%% file: rsode2.tex
\documentclass [12pt]{amsart}

\usepackage{amssymb,amsxtra,amsfonts}
\usepackage{epsf}              %
\usepackage{epsfig}             %
\usepackage{graphics}

\def\part#1{\frac{\partial\phantom{q}}{\partial#1}}

\makeatletter
 \newlength{\typesize}
\makeatother

\newlength{\vvoff}
\newlength{\hhoff}

\newcommand\beq{\begin{equation}}
\newcommand\eeq{\end{equation}}

\newcommand{\IG}{\mathbb{G}}

\newcommand{\IN}{\mathbb{N}}
\newcommand{\IP}{\mathbb{P}}

\newcommand{\IC}{\mathbb{C}}
\newcommand{\IZ}{\mathbb{Z}}

\newcommand{\IT}{\mathbb{T}}
\newcommand{\IV}{\mathbb{V}}

\newcommand{\cM}{\mathcal{M}}

\newcommand{\cD}{\mathcal{D}}

\newcommand{\cN}{\mathcal{N}}
\newcommand{\cO}{\mathcal{O}}
\newcommand{\cQ}{\mathcal{Q}}

\newcommand{\g}{       \mathfrak{g}     }

\newcommand{\lb}{       \mathfrak{b}            }

\newcommand{\gl}{       \mathfrak{gl}     } %

\newcommand{\id}{\text{\rm Id}} 
\newcommand{\Id}{\text{\rm Id}}
\newcommand{\union}{\cup}

\newcommand{\h}{\mathfrak{h}}

\newcommand{\bd}{{\bf d}}

\newcommand{\pf}{\begin{bpf}}

\newcommand{\pfms}{\begin{bpfms}}
\newcommand{\epf}{\end{bpf}\hfill$\square$\\}           %
\newcommand{\epfms}{\end{bpfms}\hfill$\square$\\}               %

\newcommand{\idea}{\begin{bidea}}

\newcommand{\eidea}{\end{bidea}\hfill$\square$\\}           %

\newcommand{\sk}{\begin{bsk}}    %

\newcommand{\esk}{\end{bsk}\hfill$\square$\\}           %
\newcommand{\sketch}{\begin{bsketch}}%

\newcommand{\esketch}{\end{bsketch}\hfill$\square$\\}

\newcommand{\wt}{\widetilde}

\newcommand{\wh}{\widehat}
\newcommand{\al}{\alpha}

\newcommand{\be}{\beta}

\newcommand{\la}{\lambda}
\newcommand{\La}{\Lambda}
\newcommand{\si}{\sigma}

\newcommand{\Lie}{\mathop{\rm Lie}}

\newcommand{\ad}{{\mathop{\rm ad}}}
\newcommand{\rank}{\mathop{\rm rank}}

\newcommand{\tr}{{\mathop{\rm Tr}}}

\newcommand{\Hom}{\mathop{\rm Hom}}

\newcommand{\PSL}{{\mathop{\rm PSL}}}
\newcommand{\GL}{{\mathop{\rm GL}}}

\newcommand{\Image}{\mathop{\rm Im}}

\newcommand{\End}{\mathop{\rm End}}

\newcommand{\Stab}{{\mathop{\rm Stab}}}

\newcommand {\eps}{\varepsilon}

\newcommand{\spq}{/\!\!/}

\def\mapright#1{\smash{
        \mathop{\longrightarrow}\limits^{#1}}}

\def\hookmapright#1{\smash{
        \mathop{\hookrightarrow}\limits^{#1}}}

\def\underset#1#2{\ \smash{\mathop{ #2 }\limits_{#1}}\ }

\newtheorem {hypo}{\bf\hspace{-\parindent}Hypothesis}
\newtheorem {thm}{Theorem}
\newtheorem {prop}[hypo]{Proposition}%
\newtheorem {lem}[hypo]{Lemma}%

\newenvironment{rmk}
        {{\bf\hspace{-\parindent}Remark.}}
        {}

\newenvironment{eg}
        {{\bf\hspace{-\parindent}Example.}}
        {}

\begin{document}
\date{}

\title[Irregular connections and Kac--Moody root systems]{Irregular Connections and\\ Kac--Moody root systems}
\author{Philip Boalch}

\maketitle

\begin{abstract}
Some moduli spaces of irregular connections on the trivial bundle over the  Riemann sphere will be identified with Nakajima quiver varieties. 
In particular this enables us to associate a Kac--Moody root system to such connections (yielding many isomorphisms between such moduli spaces, via the reflection functors for the corresponding Weyl group).
The possibility of `reading' a quiver in different ways also yields numerous isomorphisms between such moduli spaces, often between spaces of connections on different rank bundles and with different polar divisors.
Finally some results of Crawley-Boevey on the existence of stable connections will be extended to this more general context.
\end{abstract}

\section{Introduction}

It has been known since the work of Okamoto 
(see \cite{Okamoto-dynkin}) that the fourth Painlev\'e equation
has the affine Weyl group of type $\wh A_2$ as symmetry group.

\begin{equation}
\begin{matrix}
	\input{affineA2.pstex_t}
\end{matrix}
\end{equation}\label{triangle}
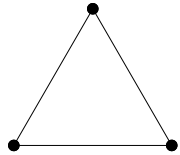
$$\text{Affine $A_2$ Dynkin diagram.}$$

Moreover the fourth Painlev\'e equation
is known (see \cite{JM81}) to be the nonlinear differential equation controlling the monodromy preserving\footnote{Tacitly this means the generalized monodromy data, including the Stokes data.} 
deformations of meromorphic connections on the trivial rank two vector bundle over the Riemann sphere of the form 
\beq\label{eq: p4conn}
\left(\frac{A}{z^3}+\frac{B}{z^2}+\frac{C}{z}\right)dz
\eeq
with $A,B,C\in \gl_2(\IC)$, and $A$ regular semisimple.
The primary aim of this article is to clarify this relationship between graphs \eqref{triangle} and connections \eqref{eq: p4conn}, 
and to extend it.
In particular this will associate a Kac--Moody root system to certain connections; recall, e.g. from \cite{kac-graphsI}, that one may associate a Kac-Moody root system to any graph with no loops connecting a vertex to itself\footnote{The Cartan matrix is $2-A$ where $A$ is the adjacency matrix of the graph.}.
On one hand this `graphical' way of viewing connections enables one to {\em see} certain isomorphisms between moduli spaces of connections which are not a priori transparent. On the other hand the Weyl group of the corresponding root systems acts (by reflection functors \cite{CB-H, Nakaj-quiver.duke, Nakaj-refl}) to yield many more isomorphisms between moduli spaces of connections.

In general complex symplectic moduli spaces of such connections on curves are obtained by fixing the formal isomorphism class at each pole (cf. \cite{smid}).
The starting point for this article was the observation
(\cite{quad} Exercise 3) that the Nakajima quiver varieties associated to the affine $A_2$ Dynkin diagram (with standard dimension vector $(1,1,1)$, i.e. the $A_2$ ALE spaces) are isomorphic to such spaces of connections
\eqref{eq: p4conn}, and the quiver  parameters match up with the choices of residues of the formal types. (We will set the real `parabolic' parameters to zero here to simplify the presentation.)
The Okamoto affine Weyl group of symmetries of the fourth Painlev\'e equation may then be understood in terms of reflection functors for the quiver variety (generating the finite Weyl group) and 
elementary/Schlesinger transforms, generating the lattice subgroup.
In these terms our primary aim is to identify some larger moduli spaces of connections with Nakajima quiver varieties for more general graphs.
This is also a generalisation of the fact, used in 
\cite{CB-additiveDS}, that moduli spaces of Fuchsian systems on the Riemann sphere are quiver varieties for star-shaped quivers.

In rough terms 
(a weak version of) the main result of this article is as follows.
Let $k\ge 1$ be an integer and let $I_1,\ldots,I_k$ be finite sets.
Let $\Gamma$ be the complete $k$-partite graph on the sets $I_i$ 
(i.e. with nodes $I=\bigcup I_i$ and each node $e\in I_i$ is connected by 
a single edge to $f\in I_j$ if and only if $i\ne j$).
Let $\cQ$ be a quiver with underlying graph $\Gamma$.
Recall that Nakajima \cite{Nakaj-quiver.duke} has defined quiver varieties $\cN_\cQ(\bd, \la)$ 
attached to any quiver (given the choice of a dimension vector 
$\bd\in \IZ_{\ge 0}^I$ and some parameters 
$\la\in \IC^I$ with $\la\cdot \bd=0$), which are smooth hyperk\"ahler manifolds for generic parameters.

\begin{thm}\label{thm: intro}
For any such quiver $\cQ$ 
there is a moduli space $\cM^*$  of meromorphic connections with fixed formal types on a trivial holomorphic vector bundle over the Riemann sphere isomorphic to a Nakajima quiver variety for the quiver $\cQ$, for some 
$\bd, \la$. 
Moreover the same holds if we first attach a leg to each node of $\,\cQ$. (A leg is a type $A$ Dynkin diagram.)
\end{thm}

In general each such quiver arises from connections in many ways, and we will explain how to `read' the connection from the quiver (and conversely how to obtain $\Gamma, \bd, \la$ from a connection).
The graphs $\Gamma$ in Theorem \ref{thm: intro}
are basically enumerated by Young diagrams, i.e. by the partition $I=\bigcup I_i$ (ignoring the legs for the moment).
One sees immediately this is a rich class of graphs, simple examples being the tetrahedral and octahedral graphs (see Figure \ref{fig: graph table}).

\ 

\begin{figure}[h]
	\centering
	\input{partite3.pstex_t}
	\caption{Graphs from partitions of $N\le 6$}\label{fig: graph table}
\center{\ (omitting the stars $\Gamma(n,1)$ and the totally disconnected graphs $\Gamma(n)$)}\label{fig: partite}
\end{figure}

Note that it is not true that all such moduli spaces of connections are quiver varieties, even on the Riemann sphere. It seems one can have at most one irregular singularity. (The simplest counterexample is the case of two poles of order two on a rank two bundle; then  $\cM^*$ is a 
$D_2$ ALF space.)
We will show however that all such moduli spaces arise from such quivers in the case of connections having one pole of order at most three (with semisimple irregular type) and arbitrarily many Fuchsian singularities.
(The case of a higher order pole, involving multiple edges, will be sketched in the appendix.)

A basic motivation underyling this work is a desire to better understand some of the hyperk\"ahler metrics of \cite{nabh} on the first wild nonabelian cohomology of a curve (the wild Hitchin moduli spaces).
In particular, on the Riemann sphere, such spaces are approximated by spaces of connections on the holomorphically trivial bundle (cf. 
\cite{%
smid}), 
and we are showing that in some cases such approximations are Nakajima quiver varieties. In turn the quiver varietes have hyperk\"ahler metrics which, since they are  hyperk\"ahler quotients of finite dimensional vector spaces, are more amenable than the metrics of
\cite{nabh}, and may be viewed as an approximation.

Physically one imagines the relation between Nakajima varieties (associated to graphs) and the Hitchin spaces (associated to Riemann surfaces) as being like degenerating a Riemann surface (string worldsheet) into a graph (Feynman diagram). In the star-shaped case this picture is presumably well-known, and we are showing it should extend to other graphs\footnote{Two intriguing features of the irregular case are: 1) that the gauge group varies at different points of the worldsheet, as one crosses the `analytic halo', and 2) that the fusion operation, involving sewing worldsheets together, appears to make sense on the {\em other side} of the `analytic halo'---(at least it does make sense on the quasi-Hamiltonian level)---as if some interactions are physically separated by the halo.}.

The possibility of `reading' the same quiver in different ways yields numerous isomorphisms between spaces of connections, often on different rank bundles.
For example taking the tetrahedral graph and fixing dimension vector $(1,2,3,4)$ yields isomorphisms between five (twelve dimensional) spaces of connections on bundles of ranks 
$6,7,8,9,10$ respectively. (In general many more isomorphisms may be obtained using reflection functors.)

These different readings generalize various known dualities, as follows.
In general a complete $k$-partite graph may be read in $k+1$ ways. 
Thus a complete bipartite graph ($k=2$) may be read in three ways, giving isomorphisms between three spaces of connections (for each set of parameters and dimension vector). Two of these three readings correspond to moduli spaces of connections having at most one pole of order $2$ (and some simple poles), and the isomorphism between these two spaces is given by Harnad duality \cite{Harn94}, which may be interpreted as a $\cD$-module Fourier--Laplace transform (cf. \cite{duality})\footnote{The third reading (as a space of connections with just one pole of order $3$) looks to be new, even in the star-shaped case.}.
In the special bipartite case when one of the sets $I_1,I_2$ contains only one element (so the graph is star-shaped with legs of length one), and if we choose a dimension vector with a $1$ at each foot, then in this case Harnad duality is the complexification of the duality that Gel$'$fand--MacPherson used (\cite{GelfMacP} p.291) to generalize the dilogarithm\footnote{It is not clear if the dilogarithm is related to the transcendental Riemann--Hilbert--Birkhoff map, passing from a connection to its monodromy/Stokes data, which is central for us.}.

By relating moduli spaces of connections to quivers one may now profit from the many algebraic results developed for quiver varieties and the deformed preprojective algebras \cite{CB-H}.
As an example in this direction we will use results of Crawley--Boevey to give precise criteria, in terms of the associated Kac--Moody root system,
 for the existence of stable connections on trivial bundles with given formal types.

To end the introduction we will mention some further directions.
There are numerous points of view on moduli spaces of connections (cf. the overview \cite{Sim96} p.11).
In particular one may take the monodromy/Stokes data of a connection 
which yields a holomorphic map $\nu$ to a moduli space $M$ of monodromy/Stokes data (the `Betti moduli space' in this context):

$$\cM^*\  \hookmapright{\nu} \ M$$

which is, at least for generic parameters, injective. 
The space $M$ may be viewed (cf. \cite{smid, saqh}) 
as a `multiplicative analogue' of the above spaces $\cM^*$ of connections.
In the cases at hand $M$ may be described directly in terms of the quiver  and one can define (\cite{salgebras}) a sort of `multiplicative quiver variety', isomorphic to $M$.
These multiplicative varieties are in general different (even for $\wh A_2, \wh A_3$) to those defined by Crawley-Boevey and Shaw \cite{CB-Shaw}, and may be viewed as a generalisation.
They may be obtained as quasi-Hamiltonian reductions of some quasi-Hamiltonian spaces associated to the quiver, related to the spaces of \cite{saqh}.
A key fact seems to be that these Betti spaces only depend on the quiver and not the particular reading.

In another direction one may now attach nonlinear differential equations to a large class of graphs (namely the isomonodromy equations for the corresponding connections, generalizing the Painlev\'e equations). 
These equations turn out to be 
essentially independent of the reading of the graph. 
This will be discussed elsewhere (\cite{gimds}), enabling us to focus on the moduli spaces here.
A related direction is to study the (linear) flat connections corresponding to these isomonodromy equations, and the quantum dualities corresponding to the different readings.
(Recall that Reshetikhin \cite{Resh92} explained how to obtain the KZ equations from the Fuchsian isomonodromy equations and \cite{bafi} Prop.4.4 explained how to obtain the extension \cite{FMTV, VTL-duke} of KZ responsible for the quantum Weyl group actions, from the simplest irregular isomonodromy equations.)

The relation between these `multi-dualities', coming from the different readings, and  
($\cD$-module) Fourier--Laplace transforms will also be described in 
\cite{gimds}. 
This relation suggests it should be possible to interpret these dualities as
Nahm transforms, as S. Szabo has done in \cite{szabo-2005} (in what we now call the bipartite case), to show that 
the metrics of \cite{nabh} 
are preserved under the dualities.

Finally, note that hyperk\"ahler moduli spaces of irregular connections on curves are crucial in Witten's approach \cite{witten-wild} to the wildly ramified case of the geometric Langlands program.

{\Small
{\bf Acknowledgements.}\  Although appearing first, this article is essentially a sequel to \cite{duality}, which studies the bipartite case in detail. I am grateful to the organizers of the Workshop on Non-linear integral transforms: Fourier-Mukai and Nahm, CRM Montreal, August 2007 for the invitation to speak, which motivated \cite{duality}.
}

{\bf Note:\ }{ Some basic/background material on quivers and Kac--Moody root systems appears in the first two appendices.}

\section{Connections on $G$-bundles}

Let $G$ be a connected complex reductive group, such as $G=\GL(V)$.
We will describe some spaces of connections on trivial holomorphic $G$-bundles over $\IP^1$. (Below a choice of invariant bilinear form on the Lie algebras of the reductive groups appearing here will be used to identify adjoint and coadjoint orbits.)

Let $\g=\Lie(G)$, choose an integer $m\ge 0$ and fix $m$ nonresonant\footnote{An element $X\in \g$ is nonresonant if $\ad_X\in\End(\g)$ has no nonzero integer eigenvalues.} adjoint orbits $\cO_i\subset \g$ and $m$ distinct points $b_i\in \IP^1\setminus\{0\}$.
Suppose $A_0\in \g$ 
is semisimple and let $H_1\subset G$ be the stabiliser of 
$A_0$ under the adjoint action.
Choose a semisimple element $A_1\in\h_1=\Lie(H_1)$
and let $H\subset H_1$ be the stabiliser of $A_1$ under the adjoint action.
Choose a nonresonant adjoint orbit $\check\cO\subset \h=\Lie(H)$.

For any $\La\in \check \cO$ define a connection

\beq \label{eq:conn}
\nabla_0=\left(\frac{A_0}{z^3}+\frac{A_1}{z^2} +\frac{\La}{z}\right)dz
\eeq

on the trivial $G$-bundle, using the inclusions $\h\subset \h_1\subset \g$.
Now consider the set $\cM^*$ 
of ($S$-equivalence classes of semistable) connections $(P,\nabla)$
where $P\to\IP^1$ is a principal $G$-bundle, isomorphic to the trivial bundle
and $\nabla$ is a meromorphic connection on $P$ with poles at $0$ and 
$\{b_i\}$ such that $\nabla$ is formally isomorphic at $0$ to $\nabla_0$ for some $\La\in \check \cO$, has simple poles at each $b_i$ with residue in $\cO_i$, and has no other poles.

In general the above moduli problem does not have a coarse moduli space which is a complex manifold. However $\cM^*$ may be obtained as a complex symplectic quotient from the following (fine) framed moduli spaces, dependent only on $A_0,A_1$, incorporating a compatible framing and allowing $\La$ to vary.
Consider the set $\wt\cM^*$ 
of isomorphism classes of connections $(P,s,\nabla)$
where 

$\bullet$
$P\to\IP^1$ is a principal $G$-bundle, isomorphic to the trivial bundle, 

$\bullet$
$\nabla$ is a meromorphic connection on $P$ formally isomorphic at $0$ to $\nabla_0$ for some nonresonant $\La\in \h$, and with simple poles at each $b_i$ with residue in $\cO_i$, and no other poles,

$\bullet$ $s\in P_0$ is a compatible framing of $(P,\nabla)$ at $0$. 

Here a framing of $P$ is simply a point of the fibre $P_0$. It is compatible if, in any local trivialisation extending $s$, $\nabla$ has the form 
$(A_0+Bz + \cdots)dz/z^3$ for some $B\in\g$ whose component in 
$\h_1=\ker(\ad_{A_0})$ equals $A_1$. The group $H$ acts transitively on the set of choices of compatible framings.

Write $G_3=G(\IC[z]/z^3)$ for the group of $3$-jets of maps from a disk to $G$ and let $B_3=\{g\in G_3 \bigl\vert\ g(0)=1\}$ be the subgroup of maps tangent to the identity.
As in \cite{smid} \S 2, the connection $\nabla_0$ determines a point of the dual of the Lie algebra $\lb_3$ of $B_3$ (which only depends on  $A_0,A_1$).

Let $\cO_B\subset \lb_3^*$ be the coadjoint orbit of $\nabla_0$.
Explicitly:
$$\cO_B = 
\left\{\left(
\frac{A_0}{z^3}+\frac{A_1+[X,A_0]}{z^2}\right)dz\ \bigl\vert \  X\in \g\ 
\right\}\cong \Image(\ad_{A_0})\subset \g.
$$

Since it is a coadjoint orbit $\cO_B$ has a complex symplectic structure (and in this case in fact has global Darboux coordinates).
Also, since  $H$ fixes $A_0,A_1$, there is a natural action of $H$ on 
$\cO_B$. As in Lemma 2.3 of \cite{smid} this action is Hamiltonian with moment map given by $-\La\in\h$ (i.e. an element of $\cO_B$ may be viewed as a connection on the disc and then transformed to a connection with principal part \eqref{eq:conn} for a unique $\Lambda\in \h$, and this defines a map 
$\La:\cO_B\to \h$).

\begin{prop}
The framed moduli space $\wt\cM^*$ is $H$-equivariantly 
isomorphic to the product
$$\cO_B\times \cO_1\times\cdots\times\cO_m$$
where $H\subset G$ acts diagonally.
\end{prop}
\pf This is similar to \cite{smid} \S2. 
\epf

Now $\cM^*$ is obtained by performing the complex symplectic quotient by $H$ at the orbit $\check \cO$ (i.e. by restricting $\La$ to be in $\check\cO$ and forgetting the framing).

\begin{rmk}
Connections on trivial bundles are just Lie algebra valued one-forms, so may equally be interpreted as Higgs fields. We prefer the connection viewpoint here though since in the sequel  \cite{salgebras} we will be considering monodromy and Stokes data. Note however that the duality of Harnad is the nonautonomous version of an earlier duality for Higgs fields (see \cite{AHH-dual}), the first case of which seems to be Moser's alternative Lax pair for the $n$-dimensional rigid body.
Similarly the `multi-dualities' of this article apply equally in the autonomous context.
\end{rmk}

\section{Relation to quivers}

Now we will
specialise to the case $G=\GL(V)$ for some complex vector space $V$, and explain the relation to quivers.
(A glance at the appendix on quiver varieties may be useful before reading this section.)

First we consider the case $m=0$ when there are no simple poles.

The data $A_0,A_1$ determine a direct sum decomposition 
$V=\bigoplus_{i\in I} V_i$ such that $H=\prod_{i\in I}\GL(V_i)\subset G$.
Namely the $V_i$ are the simultaneous eigenspaces of $A_0,A_1$.
Moreover the index set $I$ is partitioned according to the eigenspace decomposition of $A_0$. Namely if $A_0$ has $k$ distinct eigenvalues 
then $I=\bigcup_1^k I_j$ for disjoint subsets $I_j\subset I$ such that the 
$j$th eigenspace of $A_0$ is $W_j=\bigoplus_{i\in I_j}V_i\subset V$.

Now let $\Gamma_c$ be the complete $k$-partite graph on the sets $I_i$.
Thus $I$ is the set of nodes of $\Gamma_c$, and nodes $e\in I_i$ and $f\in I_j$ are connected by a single edge if and only if $i\ne j$.
Note each node corresponds to a vector space $V_i$.

Choose an ordering of $I$ and orient $\Gamma_c$ accordingly and let $\cQ_c$ be the corresponding quiver.
Consider the symplectic vector space
$$\IV(\cQ_c)=T^*\bigoplus_{i<j}\Hom(V_i,V_j)$$
consisting of maps in both directions between $V_i$ 
and $V_j$ for all $i\ne j$.
This has a natural Hamiltonian action of $H$.

\begin{lem} The coadjoint orbit $\cO_B$ is isomorphic to $\IV$ as a Hamiltonian $H$-space.
\end{lem}
\pf Direct computation. \epf

(Note that the isomorphism depends on $A_0$.)
In the terminology of the appendix, $\cO_B$ is thus an `open' quiver variety.

Now the orbit $\check \cO\subset \h$ is isomorphic to a product of orbits
$\check\cO_i\subset \End(V_i)$ (since $H=\prod\GL(V_i)$).
Each orbit $\check\cO_i$ corresponds to a leg of a quiver, possibly consisting of just one node (see \S\ref{subsection: legs}), with one open node.
In this correspondence each node of the leg is assigned a complex number and a finite dimensional vector space (with $V_i$ appearing at one end of the $i$th leg).

Thus we may construct a larger quiver $\cQ$ by gluing the open node of the $i$th such leg to the $i$th node of $\Gamma_c$.
The following is now straightforward.

\begin{prop} If all the orbits $\check\cO_i$ are closed then the moduli space $\cM^*$ (the symplectic quotient of $\cO_B$ by $H$ at $\check\cO$) is isomorphic to the quiver variety of $\cQ$ with vector spaces and parameters as assigned on each leg.
\end{prop}

This establishes the result of the introduction (in particular it is sufficient to consider connections just with one pole of order three, and no others on $\IP^1$).

Now consider the case $m>0$, so we have some adjoint orbits $\cO_i\subset 
\End(V)$ for $i=1,\ldots,m$.
Each of these corresponds to a leg, but we wish to consider the action of 
$H\subset G$ on $\cO_i$ rather than all of $G$ and this can be arranged by splaying the end of the leg according to $V=\bigoplus V_i$ 
(see \S\ref{subsection: legs}).
Thus we have $m$ splayed legs and we glue these on to $\Gamma_c$ in the obvious way (matching up the vector spaces $V_i$).

This is arranged so that if we do the symplectic quotient defining the quiver variety at all the nodes of the legs, except those now on $\Gamma_c$, the result is
just the product 
$$\cO_B\times \overline\cO_1\times\cdots\times\overline\cO_m$$
where the bar denotes the taking the closure of the orbit.
The group $H$ remains as the automorphisms of the vector spaces at the remaining nodes (those of $\Gamma_c$), and we can proceed as above and perform the symplectic quotient by $H$ at $\check\cO$, to obtain $\cM^*$  (if all the orbits are closed).

Notice that the class of graphs obtained in this way is no larger than in the case $m=0$;
The graphs obtained for $m=0$ are determined by the partition of $I$ and the lengths of the (unsplayed) legs we glue on (we will refer to these are `simple' legs).
Write $N=\# I$, which is partitioned in to $k$ parts of sizes $n_i=\# I_i$. Let $$\Gamma(n_1,\ldots,n_k)$$ denote the corresponding complete $k$-partite graph, with $N=\sum n_i$ nodes. 

A basic observation is that  gluing $m$ splayed legs 
on to the graph    
$\Gamma(n_1,\ldots,n_k)$
is the same as gluing $m$ simple legs (each of length one less) on to the 
$m$ nodes in the first `part' of the larger graph
$$m\cdot \Gamma(n_1,\ldots,n_k) := \Gamma(m, n_1,\ldots,n_k).$$
This observation enables us to obtain many isomorphisms between
moduli spaces.

Suppose $\Gamma$ is a graph obtained by gluing $N$ simple legs (of lengths $l_i \ge 0$, $i=1,\ldots,N$) on to a complete $k$-partite graph $\Gamma(n_1,\ldots,n_k)$.
Then we will refer to $\Gamma_c=\Gamma(n_1,\ldots,n_k)$ as the {\em centre} of the graph. Thus the centre has $N=\sum n_i$ nodes and $\Gamma$ has $\sum l_i$ nodes. (The centre is uniquely determined except in the star-shaped case.)

\begin{eg}
The star-shaped case occurs when the centre of the graph is bipartite and one of the parts in the partition has just one element (this element corresponds to the central node), i.e. the centre of the graph is $\Gamma(n,1)$ for some $n$.
In this case there is some ambiguity\footnote{E.g.
if $\Gamma$ is a star with $3$ legs of length $2$ (affine $E_6$) the
centre is either $\Gamma(1,3)$ or $\Gamma(1,2)$, and $\Gamma$ arises respectively by adding either $3$ simple legs of length $1$, or $2$ such legs plus a leg of length $2$ to the central node.} 
 in choosing the centre of the quiver, but once this choice is made the three possible readings are as in the following example.
\end{eg}

\begin{eg}
The general bipartite case $k=2$ has three readings, and is special since in two readings one may take $A_0=0$ so the connection has only a pole of order $2$, as follows. The centre of the graph has the form $\Gamma(p,q)$ and in  general $\Gamma(p,q)$ has  three readings:
$$\Gamma(p,q)= p\cdot\Gamma(q) = q\cdot\Gamma(p).$$
The first case corresponds to connections having a pole of order three and no others on $\IP^1$, and the leading term $A_0$ should have two distinct eigenvalues (and $A_1$ should have exactly $p$, respectively $q$, distinct eigenvalues in the two eigenspaces of $A_0$, with arbitrary multiplicities).
In the other two cases $A_0$ should be any scalar matrix (one eigenvalue)---thus in particular we may take $A_0=0$ so the corresponding connection has a pole of order just two, with leading coefficient $A_1$, and $A_1$  should have $q$ (resp. $p$) distinct eigenvalues. These connections should also have $p$  (resp. $q$) Fuchsian singularities. (The star-shaped case is special since then one of $p$ or $q$ equals $1$, and so we may take $A_1=0$ in one reading, to obtain a Fuchsian system.)
This `bipartite' quiver description of such connections having just one pole of order $2$ 
is basically equivalent to that of Jimbo--Miwa--Mori--Sato \cite{JMMS} (Appendix $5$) and Harnad \cite{Harn94}  (the novelty here being the attaching of the legs to fix the coadjoint orbits).  This description was used by Harnad to define a duality between such connections which we now interpret, as the above two (latter) readings of the bipartite quiver.
This is discussed further in \cite{duality}.
(A key observation of \cite{duality} is that since the square is a complete bipartite graph, in effect \cite{JMMS, Harn94} are telling us how to construct connections {\em directly} from the affine $A_3$ Dynkin diagram.)\footnote{The first reading (just one pole of order $3$) seems to be new even in the simplest examples: E.g. for the square $(p,q)=(2,2)$ with dimension vector $(1,1,1,1)$ this yields a rank four Lax pair for the fifth Painlev\'e equation, and similarly, but with rank six, for the sixth Painlev\'e equation, with the star $(p,q)=(4,1)$ and dimension vector $(1,1,1,1,2)$, the $2$ being at the central node.}
\end{eg}

\begin{eg}
Consider the tetrahedral graph $\Gamma(1,1,1,1)$. This may be viewed as gluing a splayed leg of length one ($\Gamma(3,1)$) on to any face (a triangle $\Gamma(1,1,1)$).
Thus for example 
if the vector spaces at the nodes have dimensions $1,2,3,4$,
removing each of the nodes in turn realizes the corresponding quiver varieties as spaces of connections with $1$ simple pole and a pole of order $3$ on bundles of rank $6,7,8,9$ respectively (the rank is the sum of the dimensions of the vector spaces on the remaining face). On the other hand if we remove no nodes, it is a space of connections with just one pole of order $3$ on a rank $10$ bundle. 
Alternatively if we take dimension vector $1,1,1,1$, removing different nodes yields isomorphisms between spaces of connections of the same type (pole orders $1$ and $3$ on bundles of rank $3$).
\end{eg}

\begin{eg} Similarly the triangle $\Gamma(111)$ may be read in various ways. In the standard reading $1\cdot\Gamma(11)$ one obtains connections with one pole of order $3$ and one pole of order $1$. (If the dimension vector is $(1,1,1)$ this is the space of connections in the introduction, on a rank $2$ bundle, related to the fourth Painlev\'e equation.) The alternative reading is as a space of connections with a pole of order $3$ and no others. 
(In the case of dimension vector $(1,1,1)$ these connections are now on a bundle of rank $3=1+1+1$; they live in a moduli space of complex dimension two, and this alternative reading appears in an explicit form in \cite{JKT07}).
\end{eg}

\begin{figure}[h]
	\centering
	\input{readings221.pstex_t}
	\caption{Four readings of 
$\Gamma(221)=2\cdot\Gamma(21)=1\cdot\Gamma(22)$}\label{fig: 221readings}
\end{figure}
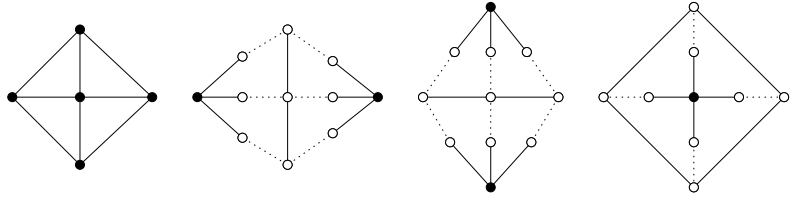

\begin{eg}
Similarly the graph $\Gamma(221)$ may be read in four ways, as indicated in Figure \ref{fig: 221readings}, corresponding to connections with poles of orders $3, 3+1+1, 3+1+1$ and  $3+1$ respectively.
\end{eg}

\ 

\begin{eg}
Consider the graph $\Gamma(3,2,1)$. In the first instance this describes a space of connections having just one pole of order $3$. Then removing each part in turn yields isomorphisms between this space and spaces of connections with a pole of order $3$ and respectively $3, 2$ or $1$ simple poles.
\end{eg}

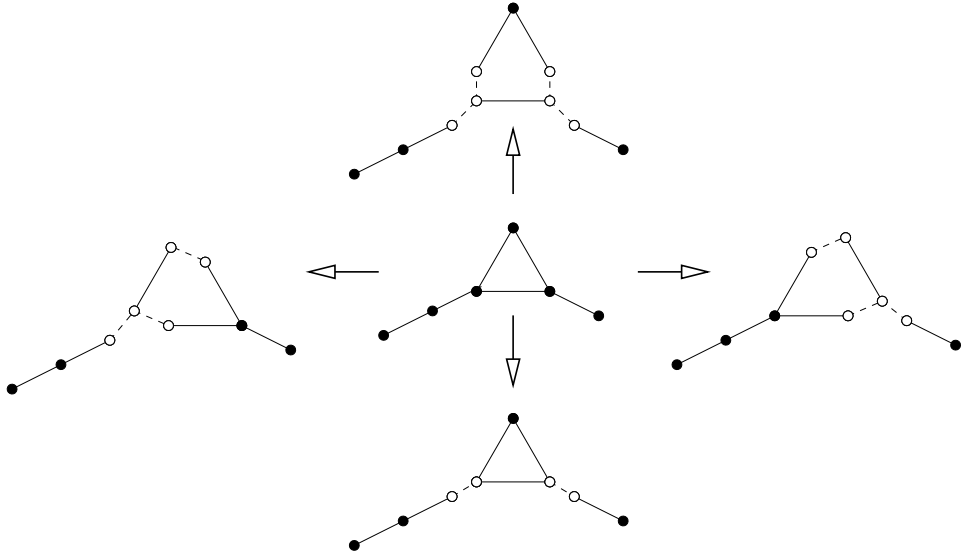
\begin{figure}[h]
	\centering
	\input{readings-legs3.pstex_t}
	\caption{Deconstructing a triangle with legs.}\label{fig: readings-legs}
\end{figure}

\begin{eg}
The examples so far have ignored the possibility of gluing some simple legs to each node of the centre of the graph. This is illustrated in Figure \ref{fig: readings-legs}, where the central piece is a triangle 
($\Gamma(111)$) in the reading at the bottom and an interval ($\Gamma(11)$) in the other three readings.
\end{eg}

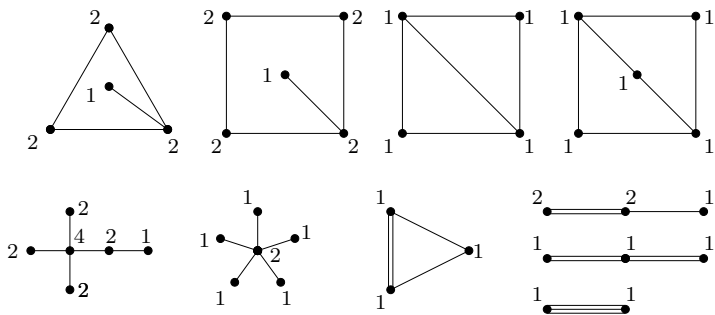
\begin{figure}[h]
	\centering
	\input{dim4.pstex_t}
	\caption{Some four dimensional cases.}\label{fig: dim4}
\end{figure}

\begin{rmk}
Recall that the `next simplest' (and most studied) class of Kac--Moody algebras after the affine algebras are the hyperbolic ones, i.e. those such that any proper subdiagram of their Dynkin diagram is either an affine or finite Dynkin diagram.
(In particular there is a simpler characterisation of their roots.) 
Thus these should correspond to the next simplest moduli spaces of connections after those of complex dimension two, which are related to affine algebras.
For example finding spaces of stable connections of complex dimension $4$ then (after Theorem \ref{thm: excons} below) basically amounts to finding integral vectors of norm $-2$.
(There are infinitely many hyperbolic diagrams that arise in the context of the present paper cf. \cite{wan}; 
five corresponding to the graphs $\Gamma(1111), \Gamma(211), \Gamma(32)$ and the two graphs obtained by attaching a single leg of length one to the square or the triangle, plus
$5$ star-shaped diagram, $5$ with double bonds---see the appendix---and an infinite family with just two nodes and a single higher order edge.)
For example one may always take an affine ADE Dynkin diagram with dimension vector the minimal imaginary root $\delta$, then double $\delta$ and glue a single leg of length one (with dimension one at the foot) on to the extending node, to obtain a diagram with a dimension vector for a quiver variety of dimension $4$. There are other examples however, see Figure \ref{fig: dim4}.

\end{rmk}

\section{Connections and quivers}

In the section we will give a more formal presentation.
(In particular we will pay more attention to orientations and signs---in reality one needs to adjust the orientation to read a quiver in different ways, and this introduces some signs---these are necessary to preserve the complex symplectic structures.)

\subsection{Abstract data $(I,J)$.}

Let $J$ be a nonempty, totally ordered, finite set and let $I_j$ be nonempty finite sets for each $j\in J$. 
Let $I_0$ be another (possibly empty) finite set.
Let $$I = I_0\cup \bigcup_{j\in J} I_j$$
be the disjoint union of the $I_j$.
Thus we have a finite set $I$ and a partition of a subset of it, into parts labelled by $J$.
We will say that the data $(I,J)$ is {\em complete} if $I_0$ is empty, i.e. 
 $I = \bigcup_{j\in J} I_j$.
Otherwise it is {\em incomplete}.

\subsection{Representations of $(I,J)$.}

Define a {\em representation} of such data $(I,J)$ to be a $J$-graded finite dimensional vector space $V$---i.e. so
there is a given direct sum decomposition
$$V =\bigoplus_{j\in J} W_j,$$
together with a $I_j$-grading of each piece:
$$W_j = \bigoplus_{i\in I_j}V_i$$
and linear maps
$$\al_{ij}\in \Hom(W_i,W_j)$$
for all distinct $i,j\in J$
and elements 
$$B_i\in \End(V)$$
for all $i\in I_0$.

This data is a convenient intermediate step between connections and quivers. 
Below  we will explain how to go from here to connections, describe some basic operations on $(I,J)$ representations, and then describe how to go from here to quivers.

\subsection{Realizations of $(I,J, V, \al, {\bf B})$.}

Suppose we have a representation of $(I,J)$ as above.
Choose complex numbers $a_j\in \IC$ for $j\in J$ and
$b_i\in \IC$ for all $i\in I$.
These should be such that the $a_j$ are distinct, and such that
if $i,i'\in I_j$ for some $j\in \{0\}\cup J$ then $b_i\ne b_{i'}$.
Such a choice of complex numbers will be referred to as a choice of {\em realization data}.

Let $\Id_j\in \End(V)$ be the idempotent\footnote{i.e. the element acting as the identity on $W_j$ and as zero on the other components of the direct sum decomposition.} for $W_j\subset V$, and,
if $i\in I_j$, let $\Id_i\in \End(W_j)$ be the idempotent for $V_i$.
(We will also view $\Id_i\in \End(V)$ in the obvious way.)

Then we may define semisimple endomorphisms
$$A_0 = \sum_{j\in J} a_j\Id_j\in \End(V),$$
$$A^{(j)} = \sum_{i\in I_j}b_i\Id_i\in \End(W_j)$$
$$A_1 = \sum_{j\in J} A^{(j)}\in \End(V).$$ 

Finally define an element 
$B\in \End(V)\cong \bigoplus_{i,j\in J} \Hom(W_i,W_j)$
to have components $B_{ij}\in \Hom(W_j,W_i)$ as follows:
$B_{ij} = -\al_{ij}$ if $i>j$,
$B_{ij} = (a_j-a_i)\al_{ij}$ if $i<j$,
$B_{ij} = A^{(j)}$ if $i=j$.
Note that this uses the given total ordering of $J$.

Then, given a coordinate $w$ on the complex plane, the corresponding realization is the meromorphic connection
\beq
-\left(A_0 w + B -\sum_{i\in I_0} \frac{B_i}{w-b_i}\right) dw
\eeq
on the trivial bundle $V\times\IP^1\to \IP^1$.
This has irregular type $(A_0/z^3+A_1/z^2)dz$ at $w=\infty$, where $z=1/w$, and has residue $B_i$ at $w=b_i$ for $i\in I_0$.
Conversely given such a connection one may go in reverse to obtain data 
$I,J,V,\al, B_i$.

Note that:

1) if $(I,J)$ is complete then the realization is just 
$-\left(A_0 w + B\right) dw$, so only has a pole at $w=\infty$,

2) if $\#J=1$ then we may choose $a_j=0$ so that $A_0=0$, and the pole at $\infty$ has order at most $2$, and leading term $-A_1dw$,

3) if $\# J=1$ {\em and} $\# I_j=1$, then we may choose $A_0=A_1=0$ so the realization is a Fuchsian system.

\subsection{Passage to complete data}
Suppose we have a representation $(I,J, V, \al, {\bf B})$ of data $(I,J)$.
We will explain how to obtain a representation of some {\em complete} data $(I,\wh J)$.
This depends on some {\em twisting data} consisting of complex numbers $x_i\in \IC$ for each $i\in I_0$. We will say the twisting data is {\em proper} if $x_i$ is an eigenvalue of $B_i$ for each $i\in I_0$.

If $I_0$ is empty there is nothing to do, the data is already complete.
Otherwise define $\wh J = \{0\} \union J$, (with the ordering so that $0$ is minimal). 
For each $i\in I_0$  define a vector space
$V_i = \Image(B_i-x_i)\subset V$.
Let $q_i:V_i\to V$ be the inclusion and let $p_i=B_i-x_i:V\to V_i$.
Define 
$$W_0 = \bigoplus_{i\in I_0} V_i$$
to be the (external) direct sum.
Thus the $p_i,q_i$ are components of maps $P:V\to W_0, Q:W_0\to V$ respectively.
In turn, since $V= \bigoplus_{J} W_j$ the maps $P,Q$ have components
$$P_j:W_j\to W_0,\qquad Q_j:W_0\to W_j.$$

Thus we may define $\wh\al_{ij}\in\Hom(W_j,W_i)$ for all distinct 
$i,j\in \wh J$, namely $\wh \al_{ij}=\al_{ij}$ if $i,j\in J$ and
$$\wh\al_{0j}=P_j,\qquad \wh\al_{j0} = Q_j.$$
In other words we define
$$\wh\al = \left(\begin{matrix} 0 & P \\ Q & \al \end{matrix}\right)
\in \End(W_0\oplus V).$$
Thus if we define $\wh V = \bigoplus_{\wh J} W_j= W_0\oplus V$  we obtain a representation $(I,\wh J, \wh V, \wh \al, 0)$ of the complete data $(I,\wh J)$.

Note that $B_i = q_i\circ p_i +x_i \in \End(V)$.
To symmetrize the notation we will define: 

$$V^{(j)} = \bigoplus_{i\in\wh J\setminus\{j\}} W_i$$
for any $j\in \wh J$, so $V = V^{(0)}$ and for any $j\in \wh J$ we have
$\wh V = W_j\oplus V^{(j)}$.
In turn $\wh \al$ determines, as above, maps 
$V_i \mapright{q_i} V^{(j)}\mapright{p_i}V_i$, 
for $i\in I_j$ for any $j\in\wh J$. 

Apart from the ordering of $J$, the situation is also not perfectly symmetric since we have
\beq
p_{i} \text{\ is surjective and\ } q_{i}
\text{\ is injective} 
\eeq
for $i\in I_0$ (but not necessarily for other $i\in I$).
We will see this asymmetry disappears for stable connections/representations.

\subsection{Cycling}

Now suppose we have a representation 
$(I, \wh J, \wh V,\wh \al)$ of some complete data $(I,\wh J)$.
We wish to `cycle' the order of the elements of $\wh J$, adjusting 
$\wh \al$ as we go. Write $\wh J=(0,1,2,\ldots,l)$, 
so $0$ is the minimal element of $\wh J$ etc.
Thus $\wh V = W_0 \oplus V$ (where $V= V^{(0)}$), and 
$$\wh \al = \left(\begin{matrix} 0 & P \\ Q & \al \end{matrix}\right)
\in \End(W_0\oplus V).$$

We will replace $\wh J$ by $\wh J' = (1,2,\ldots,l,0)$ and replace $\wh \al$
by
$$\wh \al' = \left(\begin{matrix} \al & -Q \\ P & 0 \end{matrix}\right)
\in \End(V\oplus W_0) = \End(\wh V)$$
(introducing a sign so the corresponding isomorphisms of moduli spaces are symplectic).
Said differently we negate all the components $\wh \al_{j0}$ of $\wh \al$ to obtain $\wh \al'$. 
This yields a new representation $(I, \wh J', \wh V,\wh \al')$
which we will call the `cycle' of the original representation.
(This is  a cyclic operation of order $2l+2$.)

\subsection{Passage to incomplete data}

A (left) inverse to the operation of passing to complete data may be defined in the obvious way as follows.

Let $(I,\wh J,\wh V, \wh \al)$ 
be a representation of some complete data $(I,\wh J)$. 
Let $0$ be the minimal element of $\wh J$, so $\wh V = W_0\oplus V$.
Choose twisting data $x_i\in\IC$ for all $i\in I_0$.
Set $J = \wh J\setminus\{0\}$
and write 
$$\wh \al = \left(\begin{matrix} 0 & P \\ Q & \al \end{matrix}\right)
\in \End(W_0\oplus V),$$
defining $\al \in \End(V)$.
Define $B_i = q_i\circ p_i + x_i \in \End(V)$, where $p_i,q_i$ are the components of $P,Q$ respectively, as before.

This determines an incomplete representation $(I,J,V,\al,{\bf B})$, and this operation is a left inverse to the passage to complete data, provided we choose the same twisting parameters. 
If different parameters are chosen, this corresponds to adding some scalars to the residues: $B_i\mapsto B_i + y_i$.
(In terms of connections this amounts to tensoring a realization with a logarithmic connection on the trivial line bundle.)

More interestingly if we cycle and then pass between complete/incomplete data we basically obtain an action of $\IC^I$---this should be thought of as an abstract form of middle convolution in the present context\footnote{In turn these operations may be used to {\em derive} the reflection functors (cf. \cite{quad}), although for the moment we are happy to use the reflection functors as a black-box.}.

\subsection{Different realizations}
Before moving on to quivers let us describe the different 
realizations obtained from cycling.

Let $(I,\wh J,\wh V, \wh \al)$ be a complete representation.
Upon choosing  realization data we obtain a connection 
on the bundle $\wh V\times \IP^1\to \IP^1$, with just one pole at $w=\infty$.
This will be referred to as the {\em principal realization}.
(Recall that $\wh V=\bigoplus_{\wh J}W_j$.)

Now suppose we perform the cycling operation $r$ 
times $r=0,1,2,\ldots$,
to obtain $(I, \wh J^{(r)}, \wh V,\wh \al^{(r)})$.
(We will write $\wh J=(0,1,\ldots,l)$.)
If we take the  corresponding incomplete data of this
(with twisting parameters $x_i=0$ say), then
upon choosing some
realization data (not necessarily the same as that used above) 
we obtain a connection of the form
\beq
-\left(A_0 w + B - Q(w-A^{(r)})^{-1}P \right) dw
\eeq
with $\#I_r$ simple poles, on the trivial bundle $V^{(r)}\times\IP^1\to \IP^1$, where
$V^{(r)} = \bigoplus_{j\ne r} W_j$ is a  summand of $\wh V$, and
$P, Q, B$ are determined from 
$$\wh \al^{(r)} = \left(\begin{matrix} 0 & P \\ Q & \al \end{matrix}\right)
\in \End(W_r\oplus V^{(r)}) = \End(\wh V).$$
Thus we see exactly how the ranks and the number of simple poles may change.

For example, in the bipartite case we have 
$$\wh \al^{(0)} = \left(\begin{matrix} 0 & P \\ Q & 0 \end{matrix}\right),\quad
\wh \al^{(1)} = \left(\begin{matrix} 0 & -Q \\ P & 0 \end{matrix}\right)$$
and the corresponding readings are of the form
\beq
\left(-A^{(1)} + Q(w-A^{(0)})^{-1}P \right) dw
\eeq
\beq
\left(-A^{(0)} - P(w-A^{(1)})^{-1}Q \right) dw
\eeq
on the trivial bundles with fibre $W_1, W_0$ respectively, if we take $A_0=0$ and use the same realization data.
This is, up to signs\footnote{Currently (on the level of isomorphisms of the spaces $\cM^*$) we are free to adjust the realization data appropriately---this freedom will be restricted in \cite{gimds}.}, 
the duality of \cite{Harn94}.
In this bipartite case the principal realization is of the form
\beq\label{eq: prealizn}
-\left(A_0w +B\right) dw
\eeq
on the trivial bundle with fibre $\wh V=W_0\oplus W_1$, where $A_0 = a_0\id_0 + a_1 \id_1\in \End(\wh V)$ has exactly two eigenspaces and 
$B = \left(\begin{smallmatrix} A^{(0)} & cP \\ -Q & A^{(1)} 
\end{smallmatrix}\right)$ where $c=(a_1-a_0)$.
(One may show \cite{gimds} that the isomonodromy equations for such  connections of the form \eqref{eq: prealizn} are equivalent to the JMMS equations \cite{JMMS}.)
The general case should be viewed as a generalization of  this example; the key point is that the space of $(P,Q)\in T^*\Hom(W_0,W_1)$ is being generalized to the space of block off-diagonal matrices in 
$\End(\bigoplus W_j)$.

\subsection{Quivers}

Now we will obtain a representation of a quiver $\overline \cQ$ from a complete representation $(I,J,V,\al)$.
(If instead we start with an incomplete representation we may first make it complete---as above, with a choice of twisting data---and then use the procedure below; The actual quiver obtained in this way will be independent of the choice of twisting data, provided it is proper.)

Let $\Gamma_J$ be the complete graph with vertex set $J$.
Orient $\Gamma_J$, 
so that $i\to j$ if $i<j$, to obtain a quiver $\cQ_J$.  
The element $\al$ thus constitutes a representation of the double of the quiver $\cQ_J$, with the vector space $W_j$ at the $j$th node; 
$\al\in \IV(\cQ_J)$.

Now let $\Gamma_c$ be the complete $k$-partite graph on the sets $I_j$ for $j\in J$, where $k=\#J$.
This is the graph obtained by `splaying'  the $j$th node of $\Gamma_J$ by $I_j$ for each $j\in J$. 
As such it obtains an induced orientation from that of $\Gamma_J$, so is a quiver $\cQ_c$.
Moreover, putting $V_i$ at the $i$th node for all $i\in I$ we see $\IV(\cQ_c)=\IV(\cQ_J)$, so $\al$ constitutes a representation of the double of 
$\cQ_c$. 
(Specifically if $i\in I_j$ and $i'\in I_{j'}$ with $j\ne j'$ then 
the desired element of $\Hom(V_i,V_{i'})$ is $\pi_{i'}\circ q_i$ where
$\pi_{i'}:V^{(j)}\to V_{i'}$ is the projection.)

Now we will glue some (simple) legs on to $\cQ_c$.
If $i\in I_j$ define an element $\check B_i\in \End(V_i)$ 
as follows:
$$\check B_i  = p_{i}\circ \si_j \circ q_{i} $$
where $\si_j\in\End(V^{(j)})$ acts as the identity on 
$W_h$ if $h > j$ and as $-1$ on $W_h$ if $h<j$.
(This is minus the value at $\al$ of the $i$th component of the moment map for the action of $H:=\prod_{i\in I}\GL(V_i)$ on $\IV(\cQ_c)$. It also corresponds to the $i$th component of the residue of the formal type of the corresponding meromorphic connection.)
Notice that $\check B_i$ is unchanged under the cycling operation---indeed this is clear since $\check B_i$ is the $\End(V_i)$ component of 
$$-\sum_{h<j}\al_{jh}\al_{hj} + \sum_{h>j}\al_{jh}\al_{hj}\ \in \End(W_j)$$
and cycling just negates $\al_{j0}$ and reorders $J$ to move the corresponding term $\al_{j0}\al_{0j}$ over to the right-hand sum.

For $i\in I_r$, 
let $\check\cO_i\subset \End(V_i)$ be the adjoint orbit of $\check B_i$.
Let $l_i$ be one less than the degree of the minimal polynomial of $\check B_i$ and 
let $x_{ij}$ be the roots of this polynomial for $i\in I_r, j=1,\ldots, 
1+l_i$ (which we are thus choosing an order of). 
This determines a quiver leg (oriented so each arrow points towards the centre $\cQ_c$) and a representation of its double
as usual:
Define vector spaces:
$$V_{ij}= \prod_{h=1}^j(\check B_i-x_{ih})V_i,$$

let 
$q_{ij}:V_{ij}\to V_{i(j-1)}$ be the inclusion and let 
$p_{ij}:V_{i(j-1)}\to V_{ij}$  be the map given by $\check B_i-x_{ij}$, 
for $i\in I_r, j=1,\ldots, l_i$.
By construction each such $p_{ij}$ is surjective, and each $q_{ij}$ is injective, and $\check B_i=q_{i1}\circ p_{i1} + x_{i1}\in \End(V_i)$.

Now let $\cQ$ be the quiver obtained by gluing a simple leg of length $l_i$ on to the $i$th node of $\cQ_c$ for each $i\in I$.
Let $\overline \cQ$ be the double of $\cQ$.
From the data we started with we have thus constructed a representation of $\overline \cQ$, such that the maps
$p_{ij}$ are surjective (down the legs), and the maps $q_{ij}$ are injective (up the legs), for $j>0$ and for all $i\in I$.

The following moment map conditions are also satisfied.
For $i\in I$, and $j=1,\ldots,l_i$ set 
$$\la_{ij} = x_{ij}-x_{i(j+1)}.$$
Then, directly from the definition, we have 
$$q_{i(j+1)}p_{i(j+1)} - p_{ij}q_{ij} = \la_{ij}\in \End(V_{ij}).$$
Similarly at the feet, when $j=l_i$:
$$- p_{ij}q_{ij} = \la_{ij}  \in \End(V_{ij}).$$

Finally if $i\in I_r$ we have, by definition of 
$\check B_i$: 
$$q_{i1}p_{i1} - p_{i}\circ \si_r\circ q_{i} = \la_{i0} := -x_{i1}\in \End(V_i).$$
Thus in summary we have produced a point of $\mu^{-1}(\la)\subset \IV(\cQ)$ where $\la$ has components $\la_{ij}$ for $i\in I, j=0,\ldots,l_i$.

Given complete data $(I,J)$ and dimensions $d_i=\dim(V_i)$ for each $i\in I$
note that the quiver $\cQ$ is determined by the choice of orbits 
$\check \cO_i\subset \End(V_i)$, or indeed just by the lengths $l_i$.

This gives a $1-1$ correspondence between points of $\mu^{-1}(\la)$ such that the maps are injective up the legs and surjective down the legs, and
representations of (complete) 
$(I,J)$ with $\check B_i\in \check \cO_i$ for all $i$.

\begin{rmk}
Suppose we had started with some incomplete data $(I,J)$ then passed to the complete data $(I,\wh J)$ and took the quiver of that. Then the condition at the nodes $i\in I_0$ correspond to the (splayed) leg corresponding to $B_i$ (fixing its adjoint orbit), and the relation at a node $i\in I\setminus I_0$ has the form
$$\delta_i\left(\sum_{i\in I_0}(B_i-x_i) - \Lambda_0\right) = \check B_i \in 
\check \cO_i$$
where $\delta_i:\End(V)\to \End(V_i)$ is the projection, and $\Lambda_0$ is the residue of the formal type of the {\em irregular part}\footnote{i.e. the principal part minus the residue term.} at $w=\infty$ of any corresponding connection. In particular if $J=\{j\}$ has just 
one element, then $\Lambda_0=0$, and if moreover $I_j=\{i\}$ also has just one element, then 
$V=V_i$ so the relations just say that we are fixing the adjoint orbits of each $B_h$ for $h\in I_0$ as well as the orbit of their sum, i.e. this corresponds to the Fuchsian case; the additive Deligne--Simpson problem.
\end{rmk}

\subsection{Stability}

We will say that a system of ODEs (a meromorphic connection on a trivial bundle $V\times\IP^1$ over $\IP^1$) is {\em stable} if there is no invariant subsystem (i.e. if any subspace $U\subset V$ invariant by the connection is either trivial or equal to $V$).

Now a representation of data $(I,J)$ (as defined above) is equivalent to a representation of the double of the `central' 
quiver $\cQ_c$ and a  quiver representation is stable if there are no nontrivial invariant subrepresentations (recall we are setting the parabolic parameters to zero). 
It is straightforward that stable representations of $(I,J)$ correspond to stable systems.

The following result is useful when considering different `readings' of the quiver.

\begin{lem}\label{lem: inner injsurj}
In a stable representation of data $(I,J)$ all the maps $p_i$ are surjective and all the $q_i$ are injective.
\end{lem}
\pf
Without loss of generality suppose $(I,J)$ is complete.
Then the representation amounts to decompositions 
$V=\bigoplus_J W_j, W_j=\bigoplus_{I_j} V_i$, and maps 
$\al_{ij}\in \Hom(W_i,W_j)$ for all distinct $i,j\in J$.
From this one constructs maps 
$V^{(j)}\mapright{p_i} V_i\mapright{q_i} V^{(j)}$
which encode all the maps into/out of $V_i$, for all $i\in I_j$, for all $j\in J$.
Now if some $p_i$ is not surjective then we may replace $V_i$ by $\Image(p_i)$ to obtain a subrepresentation, which would be a contradiction.
If some $q_i$ is not injective, choose a nonzero subspace $K$ of its kernel,
and we get a subrepresentation by taking $K$ at node $i$ and $0$ elsewhere, which again contradicts stability.
\epf

Now the procedure described earlier extends representations of data $(I,J)$ to representations of $\overline \cQ$ such that the result has injective maps up the legs and surjective maps down the legs (and has $q_i/p_i$ injective/surjective respectively for $i\in I_0$).

\begin{lem}
Stable representations of data $(I,J)$ correspond to stable representations of $\overline \cQ$.
\end{lem}
\pf
Suppose the correspondence yields a representation of $\overline\cQ$ which has a subrepresentation. Then the dimension must drop at some node of the centre $\cQ_c$ (due to the surjectivity conditions, else it drops nowhere).
Thus it restricts to a proper subrepesentation of $\cQ_c$, i.e. of data $(I,J)$, so is zero on $\cQ_c$ and therefore everywhere. 
The converse is similar.
\epf

\subsection{Non-emptiness conditions}

Given a quiver $\cQ$  and parameters $\la$ (a complex scalar for each node), Crawley-Boevey and Holland \cite{CB-H} have defined the deformed preprojective algebra $\Pi^\la$ as the quotient by the moment map relations of the path algebra of the doubled quiver $\overline \cQ$.
For us it is important since giving a representation of $\Pi^\la$ is the same as finding a point of $\mu^{-1}(\la)$ for some dimension vector $\bf d$;
the points of the quiver variety $\cN_{\cQ}(\la,\bd)$
are in bijection with the isomorphism classes of semisimple representations of $\Pi^\la$ with dimension vector $\bd$.
Extending work of Kac,  Crawley--Boevey \cite{CB-mmap} has established precise criteria on $\bd$ in order for $\Pi^\la$ to admit a semisimple (respectively, simple) representation, i.e. for when the quiver variety is nonempty (respectively, has a stable point).
These conditions are phrased in terms of the Kac--Moody root system attached to the quiver (see the appendix for basic definitions).

\begin{thm} [\cite{CB-mmap}]\

1) $\cN_{\cQ}(\la,\bd)$ is nonempty if an only if $\bd=\sum \be_i$ is the sum of some  positive roots $\be_i$ such that $\la\cdot \be_i=0$.

2) $\cN_{\cQ}(\la,\bd)$ has a stable point if and only if $\bd$ is a positive root and moreover if for each decomposition $\bd
=\sum \be_i$ (with each $\be_i$ a positive root with $\la\cdot \be_i=0$) one has $\Delta(\bd) > \sum \Delta(\be_i)$, where $\Delta(v) := 2-(v,v)$ is the ``expected quiver variety dimension''.
\end{thm}

We may now translate this result, using our quiver interpretation of irregular connections, to give criteria for the existence of stable connections on the trivial bundle, extending
\cite{CB-additiveDS} in the Fuchsian case. 

Fix $V, I_0$ and commuting semisimple elements $A_0, A_1\in \End(V)$, and distinct points $b_i\in \IC$ for all $i\in I_0$. 
Let $H=\Stab(A_0)\cap\Stab(A_1)$.

Suppose we have a (nonresonant) 
adjoint orbit $\check \cO\subset\h\subset\End(V)$  
of $H$ and orbits $\cO_i\subset \End(V)$ for all $i\in I_0$.

These data determine, as above, a quiver $\cQ$ with dimension vector $\bd$ and parameters $\la$. We assume the orbits are chosen so that 
$\la\cdot \bd=0$, else the corresponding quiver variety will be empty, and there will be no corresponding connections---(this condition corresponds to the residue theorem for the determinant connection).

\begin{thm}\label{thm: excons}
There is a stable meromorphic connection on the bundle $V\times\IP^1$
with formal type $(A_0/z^3+A_1/z^2+\La/z)dz$ at $z=0$ with 
$\La\in\check \cO$ and 
simple poles at $w=b_i$ (where $w=1/z$) 
with residues $B_i\in\cO_i$ for all $i\in I_0$,
if and only if

1) $\bd$ is a positive root and 

2) for each decomposition $\bd
=\sum \be_i$ (with each $\be_i$ a positive root with $\la\cdot \be_i=0$) one has $\Delta(\bd) > \sum \Delta(\be_i)$, where $\Delta(v) := 2-(v,v)$.
\end{thm}

Note that the nonresonance condition is needed only to be able to phrase the result in terms of fixing the formal type---it is true in general that the given conditions ensure there is a connection whose principal part is gauge equivalent to the above form.

\pf
Most of the work has been done:
If the conditions are satisfied we have a stable point of the quiver variety and should check this corresponds to a connection of the desired form.
However the injectivity/surjectivity conditions imply the orbits are fixed as desired, as in \cite{CB-mmap}, after possibly `unsplaying' the legs for $i\in I_0$ and appealing to Lemma \ref{lem: inner injsurj}. 
The converse is clear since we have shown stable connections correspond to stable representations of $\overline \cQ$, which are simple representations of $\Pi^\la$.
\epf

In particular note that a theorem of Vinberg (see e.g. \cite{wan}) establishes the existence of positive imaginary roots for root systems which are not affine or finite.
This  may then be used to show there are {\em nonempty} moduli spaces of stable connections corresponding to all the quivers considered in this article---namely take the dimension vector to be an indivisible positive imaginary root (i.e. not an integer multiple of another positive root) and then choose $\la$ off of the finite number of hyperplanes $\la\cdot\be_i=0$ for all $\be_i$ appearing in any decomposition of the dimension vector into a sum of positive roots.

\appendix

\begin{section}{Quiver varieties}  \label{sect: quivars}

We will recall the complex symplectic approach to 
Nakajima's quiver varieties \cite{Nakaj-quiver.duke} (setting the real `parabolic' parameters to zero to simplify the notation).
It is convenient to introduce a slightly larger class of spaces allowing some `open' nodes, and to define processes of `splaying' (or `blowing up') and `gluing' of open nodes. These processes are not inverses to each other. (The open nodes are needed here only in intermediate steps and arise as moduli spaces of meromorphic connections with framings; beware they are not quite the same as the framed spaces of \cite{Nakaj-quiver.duke}.)

Let $\cQ$ be a quiver (i.e. an oriented graph) with nodes $I$ and 
edges $\cQ$. 
We will assume there are no edge loops (i.e. edges connecting a node directly to itself).

Choose a finite dimensional complex vector space $V_i$ for each $i\in I$.
Let $\bd\in \IZ^I$ be the corresponding dimension vector (i.e. the $i$th component of $\bd$ is $\dim(V_i)$).

To construct the quiver variety consider,  for each edge of $\cQ$, a homomorphism in each direction between the vector spaces at each end:

$$\IV=\IV(\cQ):=\bigoplus_{e\in \cQ} \Hom(V_{t(e)},V_{h(e)}) \oplus  \Hom(V_{h(e)},V_{t(e)})$$
where, for a directed edge $e\in \cQ$, $h(e),t(e)\in I$ are its head and tail respectively.

The vector space $\IV$ inherits a linear complex symplectic structure if we identify  $\Hom(V_{t(e)},V_{h(e)}) \oplus  \Hom(V_{h(e)},V_{t(e)})$ with the cotangent bundle  $T^*\Hom(V_{t(e)},V_{h(e)})$ in the obvious way.

Now consider the group
$$\wt \IG=\wt \IG(\cQ):= \prod_{i\in I}\GL(V_i)$$
which acts in the natural way on $\IV$. 
(If $\phi\in \Hom(V_i,V_j)$ and $g\in \wt \IG$ 
has components $g_i\in \GL(V_i)$ 
then $g\cdot \phi = g_j \circ \phi \circ g_i^{-1}$.)

The group $\wt \IG$ has a central subgroup $\IT$ isomorphic to $\IC^*$ which acts trivially on $\IV$ 
(simply take each component $g_i$ to be the same nonzero scalar).
Define $\IG = \wt \IG / \IT$ to be the quotient group (which acts effectively on $\IV$). Thus $\IG$ has a centre $Z$ of dimension $\# I-1$ (obtained by setting each $g_i$ to be an arbitrary nonzero scalar).

A moment map $\mu : \IV \to \prod_{i\in I}\gl(V_i)$ for the action of $\wt \IG$ on $\IV$ is  given by
$$\mu_i(\{\phi\}) := \sum_{e\in \cQ, h(e)=i} \phi_e \circ \phi_{e^*} -
\sum_{e\in \cQ, t(e)=i} \phi_{e^*} \circ \phi_{e}$$
where $\phi_e\in \Hom(V_{t(e)},V_{h(e)})$ and $e^*$ is the edge $e$ with reversed orientation.

Taking the sum of the traces of the components $\mu_i$ of $\mu$ clearly yields zero (since $\tr(AB)=\tr(BA)$ for any compatible rectangular matrices $A,B$) and the subspace on which this sum of traces is zero is naturally identified with the (dual of the) Lie algebra of $\IG$, and so $\mu$ is also a moment map for $\IG$.

The corresponding quiver varieties are obtained by performing the complex symplectic quotient of $\IV$ by $\IG$ at any central value 
$\lambda$ of the moment map $\mu$:
$$\cN_\cQ(\lambda, \bd) 
= \IV\underset{\la}{\spq} \IG = \mu^{-1}(\lambda)/\IG$$
where $\bd$ 
is the dimension vector.
More precisely, $\mu^{-1}(\la)$ is an affine subvariety of $\IV$ and we define $\cN_\cQ(\lambda, \bd)$ to be the affine variety associated to the ring of $\IG$ invariant functions on  $\mu^{-1}(\la)$.
(We may also view this as $\IV\spq_\la\wt \IG$ since $\IT$ acts trivially.)

The space of parameters $\la$ is given 
concretely as a subspace of 
$\prod_{i\in I}\gl(V_i)$ by the  $I$-tuples of scalar matrices whose sum of traces is zero. In other words if each diagonal entry of the $i$th matrix is 
$\lambda_i$ (so $\la\in \IC^I$) then the constraint is that
\beq\label{eq: de.la=0}
\bd\cdot\lambda := \sum_{i\in I} d_i \lambda_i = \sum_{ i\in I}\tr_{V_i}(\la_i)= 0.
\eeq

\subsection{Open nodes}

Now suppose we have a subset $I(O)\subset I$ of nodes which we will call `open'. Let 
$I(C)=I\setminus I(O)$ be the complement; the set of `closed' nodes.
Then one may repeat the above construction (with the same space $\IV$)  but only quotient by the group
$$\wt\IG(C)=\wt\IG(\cQ, C)= \prod_{i\in I(C)}\GL(V_i)\subset \wt \IG$$
of automorphisms of vector spaces at the closed nodes.
Thus we obtain spaces

$$\cN^o_\cQ(\la',\bd) = \IV \spq_{\la'} \wt\IG(C)$$

which only depend on the parameters $\la'$ on the closed nodes.
The resulting variety will still have a residual Hamiltonian 
action of the group
$$\wt\IG(O)=\wt \IG(\cQ,O)= \prod_{i\in I(O)} \GL(V_i)\subset \wt \IG.$$ 
(In other words the reduction defining $\cN_\cQ$ may be done in stages.)
Note that if all the nodes are open then the open quiver variety is just 
$\IV(\cQ)$.

\subsection{Gluing open nodes}

Suppose $\cQ$ (which may be disconnected) has two open nodes $i,j\in I(O)$ with the same dimension $d=\dim(V_{i})=\dim(V_{j})$.
Then one may glue these two nodes together to obtain 
a new quiver $\cQ'$. Label by $o$ the node of $\cQ'$ obtained from gluing $i$ and $j$, and leave it open.

Observe that the corresponding symplectic vector spaces are isomorphic:
$\IV(\cQ)=\IV(\cQ')$.
The groups acting are slightly different however.
Choose an isomorphism $V_i\cong V_j$ so it makes sense to speak of the diagonal subgroup of $\GL(V_i)\times\GL(V_j)$.
In effect we 
are then restricting  
the subgroup $\GL(V_i)\times\GL(V_j)\subset \wt \IG(\cQ)$ to its diagonal subgroup (to obtain an inclusion $\wt \IG(\cQ', O)\subset \wt \IG(\cQ, O)$).
The basic statement then is the following.

\begin{lem}
$$\cN^o_\cQ(\la,\bd) \cong \cN^o_{\cQ'}(\la,\bd')$$
as Hamiltonian $\wt \IG(\cQ', O)$-spaces.
\end{lem}
\pf
The spaces are the same since $\IV(\cQ)=\IV(\cQ')$ and the `closed' subgroups are identical.
The moment maps for the residual actions of $\wt \IG(\cQ', O)$ match up since, in full generality, the moment map for a diagonal subgroup $G\subset G\times G$ is the sum of the moment maps for each factor $G$ (the sum being the dual of the derivative of the diagonal inclusion).
\epf

(It is easy to see the above gluing process may be iterated and inverted.)

\subsection{Splaying/blowing-up nodes}

Now suppose we have a quiver $\cQ$ as above and an open node 
$o\in I(O)\subset I$.
Suppose that the vector space $V_o$ at the node $o$ has a direct sum decomposition 
$$V_o = \bigoplus_{j\in J} W_j.$$

The variety $\cN^o_\cQ(\la',\bd)$ will have an action of $\GL(V_o)\subset \wt\IG(O)$ and we wish to modify the quiver $\cQ$ such that this action (of $\GL(V_o)$) is restricted to the `block-diagonal' subgroup
$$\prod_{j\in J}\GL(W_j)\subset \GL(V_o).$$

The splaying procedure to obtain the new quiver $\cQ'$ from $\cQ$ is local at the node $o$. 
Simply replace $o$ by a node for each element of $J$ 
(so $I' = J\sqcup (I\setminus \{o\})$)
and connect up each element of $J$ in the same way as $o$ was connected to the rest of $\cQ$, i.e. repeating the edges to/from $o$ once to/from  each element of $J$.
In other words we have a surjective map
$$\pi : I'\to I$$
such that $\pi(i)=i$ if $i\in I\setminus \{o\}$ and $\pi(j)=o$ if $j\in J$.
Then the set of edges in $\cQ'$ from $i\in I'$ to $j\in I'$ is defined to be equal to the set of edges from $\pi(i)$ to $\pi(j)$ in $\cQ$.

If we define the nodes in $J$ to be open nodes of $I'$  (i.e $I'(O)=J\sqcup I(O)\setminus\{o\}$) and assign the vector space $W_j$ to each node in $J$
then (since there are no edge loops) one sees immediately that the symplectic vector spaces are the same:
$$\IV(\cQ) = \IV(\cQ').$$
However the groups acting are different, in effect we have replaced the factor $\GL(V_o)$ of $\wt\IG$ by $\prod_{j\in J} \GL(W_j)\subset \GL(V_o)$.

Some examples will appear below.
The complete $k$-partite graphs may be obtained by starting with a complete graph with $k$ nodes and then splaying each of the $k$ nodes in turn.

\subsection{Legs and orbits}\label{subsection: legs}

Let $V\cong \IC^n$ be a finite dimensional complex vector space and
let $\cO\subset \gl(V)$ be an arbitrary adjoint orbit of $n\times n$ matrices. 
We wish to recall the basic fact that the closure of $\cO$ may be obtained 
from the above construction. 
(Note that semisimple orbits are already closed.)

Suppose $A\in \cO$, let $p$ be the minimal polynomial of $A$, let $w$ be the degree of $p$ and set $l:=w-1$.
Let $x_1,\ldots,x_w\in \IC$ be the roots of $p$ (which we are thus choosing an arbitrary order of). 
We then associate to $\cO$ the linear quiver (or ``leg") with $w$ nodes, labelled $0,1,\ldots,l$ as in Figure \ref{fig: leg}.
The vertex $0$ is open and the rest are closed.

\ 

\begin{figure}[h]
	\centering
	\input{leg.pstex_t}
	\caption{A leg}\label{fig: leg}
\end{figure}

We assign dimension $n$ to the open node $0$ and the other dimensions are given by
$$d_i:= \rank(A-x_1)(A-x_2)\cdots(A-x_i).$$
Thus for example if $A$ is regular then $d_i=n-i$ for all $i$ (regardless of whether or not $A$ is semisimple).
Define the parameters to be $\la_i=x_i-x_{i+1}$ for $i=1,\ldots,l$.

\begin{lem} \label{lem: coads from legs}
The quiver variety
$\cN^o_{\cQ}(\la,\bd)$
is isomorphic to the closure of the orbit $\cO$ in $\gl(V)$. The residual action of $\GL(V)$ on this symplectic quotient corresponds to the adjoint action of $\GL(V)$ on $\cO$.
\end{lem}

\pf For nilpotent orbits this is due to Kraft--Procesi \cite{Kraft-Procesi-InvMath79}. The general case is in the appendix to 
\cite{CB-normality}. \epf

(Note also that the same result is true \cite{CB-normality} if, instead of the minimal polynomial, we take $p$ to be any polynomial with $p(A)=0$.)

Strictly speaking one does not need to consider open nodes to obtain $\cO$ but it will be convenient for us to do so, in order to obtain the action naturally.
The map associating a point of $\overline \cO\subset \gl(V)$
to a point $\phi\in\mu^{-1}(\la)\subset \IV$ (in a closed orbit of $\wt\IG(O)$) is given by $\phi\mapsto x_1 + \phi_e\circ\phi_{e^*}$, where $e$ is the edge of the leg between the nodes $1$ and $0$.
It is thus sometimes convenient to associate the parameter $\la_0=-x_1$ to the node $0$.

\subsection{Splaying a leg}

A basic example of the splaying procedure is the quiver obtained by restricting the adjoint action on an orbit $\cO$ to a block diagonal subgroup.

Suppose $V\cong \bigoplus_{j\in J} V_j$, so we have a subgroup
$\prod_{j\in J}\GL(V_j)\subset \GL(V)$.
If we splay the leg corresponding to $\cO$ according to $J$ we obtain 
the quiver as in Figure \ref{fig: splayed leg}, with $n_j=\dim(V_j)$, so $\sum n_j=n$

\ 

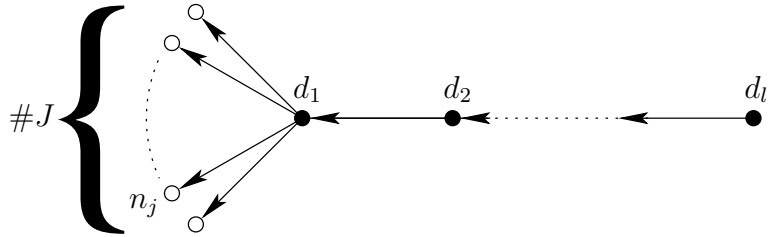
\begin{figure}[h]
	\centering
	\input{splayedleg.pstex_t}
	\caption{A leg splayed at one end}\label{fig: splayed leg}
\end{figure}

\section{Root systems and reflection functors}

\subsection{The Kac--Moody Weyl group and reflection functors}

Let $\Gamma$ be a graph with no edge loops.
One can then define a Weyl group as follows.
Let $I$ be the set of nodes and let $n=\#I$ be the number of nodes.
Define the $n\times n$ (symmetric) Cartan matrix to be 
$$C=2\ \id - A$$
where $A$ is the adjacency matrix of $\Gamma$; the $i,j$ entry of $A$ is the number of edges connecting the nodes $i$ and $j$.
The {\em root lattice} $\IZ^I=\bigoplus_{i\in I} \IZ\eps_i$ inherits a bilinear form defined by 
$$(\eps_i,\eps_j) = C_{ij}.$$

In turn we can define simple reflections $s_i$ acting on the root lattice by the 
formula
$$s_i(\be) := \be-(\be,\eps_i)\eps_i$$
for any $i\in I$.
The Weyl group is then the group generated by these simple reflections.
They satisfy (cf. \cite{kac-graphsI} p.63) the relations
$$s_i^2=1,\qquad 
s_is_j=s_js_i \text{\ if $A_{ij}=0$,}\qquad
s_is_js_i=s_js_is_j \text{\ if $A_{ij}=1$}.
$$
We can also define dual reflections $r_i$ acting on the vector space 
$\IC^I$ by the formula 
$$r_i(\la) = \la - \la_i \al_i$$
where $\la = \sum_{i\in I} \la_i\eps_i\in \IC^I$ with $\la_i\in\IC$ and where for each $i\in I$
$$\al_i := \sum_{j\in I} (\eps_i,\eps_j)\eps_j\in \IC^I.$$ 
By construction one has that $s_i(\be)\cdot r_i(\la) = \be \cdot \la$,
where the dot denotes the pairing given by 
$\eps_i\cdot\eps_j = \delta_{ij}$.
The main result we wish to quote is as follows.

\begin{thm}[\cite{Nakaj-quiver.duke, CB-H, Nakaj-refl}]
Let $\cQ$ be a  fixed quiver as above.
Then, if $\la_i\ne 0$, 
there is a natural isomorphism between the the quiver variety with dimension vector $\be$ and parameters $\la=\sum\la_i\eps_i$ 
and that with dimension vector $s_i(\be)$ and parameters $r_i(\la)$:
$$\cN_{\cQ}(\la,\be) \cong \cN_{\cQ}(r_i(\la),s_i(\be)).$$
\end{thm}

The desired reflection functors are constructed in \cite{CB-H} section $5$
(cf. \cite{CB-decomp} Lemma 2.1, \cite{CB-mmap} Lemma 2.2) and in \cite{Nakaj-quiver.duke, Nakaj-refl} in the more general 
hyperk\"ahler context. 

As an example consider the case $A_2^{++}$ of the triangle with a leg of length one attached. Then \cite{ogg84, FKN08} the index two `rotation subgroup' of the Weyl group is $\PSL_2(E)$ where $E=\IZ[\omega]$ is the ring of Eisenstein 
integers ($\omega^3=1$).

Here is a specific example. Label the nodes of $A_2^{++}$ as $1,2,3,4$ (with the $2$ in the middle, the $1$ at the foot and $2,3,4$ on the triangle).
If we take dimension vector $(1,2,2,1)$ the corresponding quiver varieties have complex dimension $2$. Indeed performing the reflections $s_1s_2s_3$ yields dimension vector $(0,1,1,1)$ so the variety is isomorphic to an $A_2$ ALE space\footnote{I am grateful to H. Nakajima for describing a similar example to me.} (which in one reading is thus also isomorphic to a space of connections on a rank $3$ bundle with $2$ poles of order $1$ and $3$).
On the other hand the Weyl group element (see  \cite{FKN08} 4.20): 
$$w=s_1s_4s_1s_2s_4s_1s_3s_1$$
has infinite order, realising the same space as a space of connections on bundles of arbitrarily high rank---indeed for $n\ge 1$ the space with dimension vector $w^n(1,2,2,1)$ may be read as a space of connections as above on bundles of rank $n^2+(n-1)+(n-2)^2$.

\begin{rmk}
Note that one also has a birational action of the group of Schlesinger transformations (i.e. the bundle modifications).
As in \cite{quad} this amounts to adding translations of the parameters $\la$ by the lattice $L=\{\la\in\IZ^I\ \bigl\vert\ \la\cdot\bd = 0 \ \}$ to obtain $W\ltimes L$.
Also as in \cite{quad} Remark 6, this becomes a biregular action upon partially compactifying to the DeRham space (allowing connections on nontrivial bundles), at least for generic parameters, when the Riemann--Hilbert--Bikhoff map is an analytic isomorphism. 
(On one hand these translations may be viewed as Hecke operators, and on another they constitute interesting nonlinear difference equations---the higher order difference Painlev\'e equations.)

Note also that in the case of an affine Dynkin diagram one may take the dimension vector to be the minimal positive imaginary root, which is fixed by all the reflections. Moreover in this case $W$ (an affine Weyl group) is  not acting effectively on the space of parameters 
$\{\la \ \bigl\vert\  \la\cdot\bd=0\}$ but acts via its quotient the finite Weyl group $W_{fin}$ (cf. \cite{quad}).
The Okamoto symmetry group is the extended affine Weyl group $W_{fin}\ltimes L$ (plus possible diagram automorphisms).
In general  $W$ itself will act effectively though, usually changing both $\bd$ and $\la$, leading to the group $W\ltimes L$ (plus diagram automorphisms).   
\end{rmk}

\subsection{The root system}
The Kac--Moody root system is the subset of the {root lattice} 
$\IZ^I=\bigoplus_{I} \IZ\eps_i$ defined as follows:
The simple roots are $\eps_i$ for $i\in I$.
The set of real roots is the Weyl group orbit of the set of simple roots.
Define the {\em fundamental region} to be the set of nonzero 
$\be\in\IN^I$ whose support is a connected subgraph of $\Gamma$ and such that $(\eps_i,\be)\le 0$ for all $i\in I$.
The set of imaginary roots is the union of the Weyl group orbit 
of the fundamental 
region and the orbit of minus the fundamental region.
The set of roots is the union of the real and imaginary roots, and a root is positive if all its coefficients are $\ge 0$.

\end{section}

\section{Higher order poles}

This appendix sketches the higher order pole case and in particular will describe the quivers that arise when one considers a pole of order 4 or more (still with semisimple irregular type).
A new feature is that edges of multiplicity greater than one occur.

Let $V$ be a finite dimensional complex vector space and
let $G=\GL(V)$.

Choose commuting semisimple elements $A_0,\ldots,A_{k-2}\in \End(V)$.
(Thus in some basis of $V$ they are all diagonal.)
Consider the irregular type:
$$\left(\frac{A_0}{z^k}+\cdots \frac{A_{k-2}}{z^2}\right)dz.$$

As usual (\cite{smid} \S2)
this may be viewed as an element of the dual of the Lie algebra of  the group $B_k$ of $k$-jets of map from a disk to $G$, tangent to the identity at the origin.
We wish to describe its  coadjoint orbit $\cO_B$ in terms of quivers.
 
Define a nested sequence of `block-diagonal' subgroups of $G$ as follows:
For $i=1,\ldots,k-1$ let
$$H_i = \Stab(A_0,\ldots,A_{i-1})=
\{g\in G\ \bigl\vert\ gA_jg^{-1}=A_j \text{\ for $j=0,\ldots,i-1$}\}.$$
Write $H_0=G$ and for convenience, $H=H_{k-1}$, so that
$$ H=H_{k-1}\subset H_{k-2}\subset \cdots \subset H_1\subset H_0=G.$$
Said differently, for each $i$, 
define  $J_i$ be the set of simultaneous eigenspaces 
in $V$ of $A_0,\ldots,A_{i-1}$ and let $V^{(i)}_j\subset V$ be the corresponding eigenspace (for $j\in J_i$).
Then, for each $i$
$$
V = \bigoplus_{j\in J_i}V^{(i)}_j,\quad 
H_i\cong \prod_{j\in J_i} \GL(V^{(i)}_j).
$$
This nesting of subspaces determines (and is determined by) a rooted tree, with nodes $J_i$ on the $i$th level, and branches according to the nesting;
$j\in J_i$ is connected to the node $l\in J_{i+1}$ if and only if
$V_l^{(i+1)}\subset V_j^{(i)}$. Said differently the branching determines (and is determined by) surjective maps
$$\pi_i: J_{i} \to J_{i-1}$$
such that  $V_j^{(i)}\subset V^{(i-1)}_{\pi_i(j)}$ for all $j\in J_i$.

Write $\h_i=\Lie(H_i)$ for the Lie algebra of $H_i$.
By construction $\h_i$ is the kernel of 
$\ad_{A_{i-1}}\vert_{\h_{i-1}}\subset \h_{i-1}$, and
there is a vector space direct sum
$$\h_{i-1}  = \h_i\oplus \h'_i$$
where $\h'_i = 
\Image(\ad_{A_{i-1}}\vert_{\h_{i-1}}) = [A_{i-1},\h_{i-1}]\subset \h_{i-1}$.
In terms of eigenspaces we have
$$\h'_i = \bigoplus_{\begin{smallmatrix}j\ne l\in J_i \\ \pi_i(j)=\pi_i(l)\end{smallmatrix}} 
\Hom(V_j^{(i)},V_l^{(i)}) \subset \h_{i-1}.$$

As a space one has $\cO_B\cong \bigoplus_{i=1}^{k-2} (\h'_i)^{k-1-i}$
(this follows from the standard procedure to put a connection germ in normal form).
The above description of $\h'_i$ enables one to encode this as an (open) quiver variety.
This is most easily described in terms of a slight generalization of the (above) notion of splaying an open node of a quiver.
The `$r$-fission' of an open node of a quiver by a set $J$ is defined as follows. First 
splay the node as usual (replace the node by a new node for each element of $J$, and connect up each node of $J$ in the same way as the original node was connected to the rest of the quiver). 
Then add new edges joining each distinct pair of nodes of $J$ by exactly $r$ edges. Thus the usual splaying process is the case $r=0$.

Thus $\cO_B$ determines a graph as follows.
Start with a single node, and perform the $(k-2)$-fission by the set $J_1$ (to obtain a complete  graph with nodes $J_1$, and each edge of multiplicity $k-2$).
Then, for each $i\in J_1$, perform $(k-3)$-fission on the node $i$ 
by the set $\pi_2^{-1}(i)\subset J_2$. The resulting graph has nodes $J_2$.
Then repeat this process until the $0$-fissions are performed to obtain a graph $\Gamma$ with nodes $J_{k-1}$, and thus a quiver $\cQ$ upon choosing an orientation. Recall each node $i\in J_{k-1}$ is associated to a vector space $V_i^{(k-1)}$.

{\bf Claim.\ } $\cO_B\cong \IV(\cQ)$ as Hamiltonian $H$-spaces.

Then spaces of connections with one pole of order $k$ and some Fuchsian singularities may  be obtained as quiver varieties as before, by gluing on a splayed leg for each Fuchsian singularity and a simple leg for each node of $J_{k-1}$.

The reader may readily construct numerous examples;
in the simplest case, of one pole of order four on a rank two bundle with $A_0$ regular semisimple, one obtains  
the affine $A_1$ Dynkin diagram, as expected from 
\cite{Okamoto-dynkin, JM81}.

\renewcommand{\baselinestretch}{1}              %
\normalsize
\bibliographystyle{amsplain}    \label{biby}
\bibliography{../thesis/syr} 

\vspace{0.5cm}   
\'Ecole Normale Sup\'erieure et CNRS, 
45 rue d'Ulm, 
75005 Paris, 
France

www.dma.ens.fr/$\sim$boalch

boalch@dma.ens.fr \qquad \qquad \qquad \qquad \qquad \quad\,\,  

\

\end{document}

%% file: affineA2.pstex_t
\begin{picture}(0,0)%
\includegraphics{affineA2.pstex}%
\end{picture}%
\setlength{\unitlength}{2901sp}%
\begingroup\makeatletter\ifx\SetFigFontNFSS\undefined%
\gdef\SetFigFontNFSS#1#2#3#4#5{%
  \reset@font\fontsize{#1}{#2pt}%
  \fontfamily{#3}\fontseries{#4}\fontshape{#5}%
  \selectfont}%
\fi\endgroup%
\begin{picture}(1456,1274)(173,-3489)
\end{picture}%

%% file: partite3.pstex_t
\begin{picture}(0,0)%
\includegraphics{partite3.pstex}%
\end{picture}%
\setlength{\unitlength}{1865sp}%
\begingroup\makeatletter\ifx\SetFigFontNFSS\undefined%
\gdef\SetFigFontNFSS#1#2#3#4#5{%
  \reset@font\fontsize{#1}{#2pt}%
  \fontfamily{#3}\fontseries{#4}\fontshape{#5}%
  \selectfont}%
\fi\endgroup%
\begin{picture}(9731,9351)(2198,-9359)
\put(2701,-4336){\makebox(0,0)[lb]{\smash{{\SetFigFontNFSS{7}{8.4}{\rmdefault}{\mddefault}{\updefault}{$32$}%
}}}}
\put(4726,-4336){\makebox(0,0)[lb]{\smash{{\SetFigFontNFSS{7}{8.4}{\rmdefault}{\mddefault}{\updefault}{$221$}%
}}}}
\put(6751,-4336){\makebox(0,0)[lb]{\smash{{\SetFigFontNFSS{7}{8.4}{\rmdefault}{\mddefault}{\updefault}{$311$}%
}}}}
\put(8776,-4336){\makebox(0,0)[lb]{\smash{{\SetFigFontNFSS{7}{8.4}{\rmdefault}{\mddefault}{\updefault}{$2111$}%
}}}}
\put(10801,-4336){\makebox(0,0)[lb]{\smash{{\SetFigFontNFSS{7}{8.4}{\rmdefault}{\mddefault}{\updefault}{$11111$}%
}}}}
\put(5851,-9286){\makebox(0,0)[lb]{\smash{{\SetFigFontNFSS{7}{8.4}{\rmdefault}{\mddefault}{\updefault}{$3111$}%
}}}}
\put(6751,-6811){\makebox(0,0)[lb]{\smash{{\SetFigFontNFSS{7}{8.4}{\rmdefault}{\mddefault}{\updefault}{$2211$}%
}}}}
\put(10576,-6811){\makebox(0,0)[lb]{\smash{{\SetFigFontNFSS{7}{8.4}{\rmdefault}{\mddefault}{\updefault}{$21111$}%
}}}}
\put(8776,-6811){\makebox(0,0)[lb]{\smash{{\SetFigFontNFSS{7}{8.4}{\rmdefault}{\mddefault}{\updefault}{$222$}%
}}}}
\put(2701,-6811){\makebox(0,0)[lb]{\smash{{\SetFigFontNFSS{7}{8.4}{\rmdefault}{\mddefault}{\updefault}{$42$}%
}}}}
\put(4726,-6811){\makebox(0,0)[lb]{\smash{{\SetFigFontNFSS{7}{8.4}{\rmdefault}{\mddefault}{\updefault}{$411$}%
}}}}
\put(3826,-9286){\makebox(0,0)[lb]{\smash{{\SetFigFontNFSS{7}{8.4}{\rmdefault}{\mddefault}{\updefault}{$321$}%
}}}}
\put(3736,-1861){\makebox(0,0)[lb]{\smash{{\SetFigFontNFSS{7}{8.4}{\rmdefault}{\mddefault}{\updefault}{$111$}%
}}}}
\put(7786,-1861){\makebox(0,0)[lb]{\smash{{\SetFigFontNFSS{7}{8.4}{\rmdefault}{\mddefault}{\updefault}{$211$}%
}}}}
\put(5761,-1861){\makebox(0,0)[lb]{\smash{{\SetFigFontNFSS{7}{8.4}{\rmdefault}{\mddefault}{\updefault}{$22$}%
}}}}
\put(9811,-1861){\makebox(0,0)[lb]{\smash{{\SetFigFontNFSS{7}{8.4}{\rmdefault}{\mddefault}{\updefault}{$1111$}%
}}}}
\put(7831,-9286){\makebox(0,0)[lb]{\smash{{\SetFigFontNFSS{7}{8.4}{\rmdefault}{\mddefault}{\updefault}{$33$}%
}}}}
\put(9676,-9286){\makebox(0,0)[lb]{\smash{{\SetFigFontNFSS{7}{8.4}{\rmdefault}{\mddefault}{\updefault}{$111111$}%
}}}}
\end{picture}%

%% file: readings221.pstex_t
\begin{picture}(0,0)%
\includegraphics{readings221.pstex}%
\end{picture}%
\setlength{\unitlength}{2486sp}%
\begingroup\makeatletter\ifx\SetFigFontNFSS\undefined%
\gdef\SetFigFontNFSS#1#2#3#4#5{%
  \reset@font\fontsize{#1}{#2pt}%
  \fontfamily{#3}\fontseries{#4}\fontshape{#5}%
  \selectfont}%
\fi\endgroup%
\begin{picture}(7801,1906)(-97,-1914)
\end{picture}%

%% file: readings-legs3.pstex_t
\begin{picture}(0,0)%
\includegraphics{readings-legs3.pstex}%
\end{picture}%
\setlength{\unitlength}{2693sp}%
\begingroup\makeatletter\ifx\SetFigFontNFSS\undefined%
\gdef\SetFigFontNFSS#1#2#3#4#5{%
  \reset@font\fontsize{#1}{#2pt}%
  \fontfamily{#3}\fontseries{#4}\fontshape{#5}%
  \selectfont}%
\fi\endgroup%
\begin{picture}(8791,5055)(-2752,-12803)
\end{picture}%

%% file: dim4.pstex_t
\begin{picture}(0,0)%
\includegraphics{dim4.pstex}%
\end{picture}%
\setlength{\unitlength}{2155sp}%
\begingroup\makeatletter\ifx\SetFigFontNFSS\undefined%
\gdef\SetFigFontNFSS#1#2#3#4#5{%
  \reset@font\fontsize{#1}{#2pt}%
  \fontfamily{#3}\fontseries{#4}\fontshape{#5}%
  \selectfont}%
\fi\endgroup%
\begin{picture}(8078,3529)(2776,-3509)
\put(3016,-1591){\makebox(0,0)[lb]{\smash{{\SetFigFontNFSS{8}{9.6}{\rmdefault}{\mddefault}{\updefault}{$2$}%
}}}}
\put(9181,-1636){\makebox(0,0)[lb]{\smash{{\SetFigFontNFSS{8}{9.6}{\rmdefault}{\mddefault}{\updefault}{$1$}%
}}}}
\put(8821,-3346){\makebox(0,0)[lb]{\smash{{\SetFigFontNFSS{8}{9.6}{\rmdefault}{\mddefault}{\updefault}{$1$}%
}}}}
\put(6166,-2626){\makebox(0,0)[lb]{\smash{{\SetFigFontNFSS{8}{9.6}{\rmdefault}{\mddefault}{\updefault}{$1$}%
}}}}
\put(4996,-2716){\makebox(0,0)[lb]{\smash{{\SetFigFontNFSS{8}{9.6}{\rmdefault}{\mddefault}{\updefault}{$1$}%
}}}}
\put(5176,-3391){\makebox(0,0)[lb]{\smash{{\SetFigFontNFSS{8}{9.6}{\rmdefault}{\mddefault}{\updefault}{$1$}%
}}}}
\put(2791,-2851){\makebox(0,0)[lb]{\smash{{\SetFigFontNFSS{8}{9.6}{\rmdefault}{\mddefault}{\updefault}{$2$}%
}}}}
\put(3601,-3301){\makebox(0,0)[lb]{\smash{{\SetFigFontNFSS{8}{9.6}{\rmdefault}{\mddefault}{\updefault}{$2$}%
}}}}
\put(3601,-3301){\makebox(0,0)[lb]{\smash{{\SetFigFontNFSS{8}{9.6}{\rmdefault}{\mddefault}{\updefault}{$2$}%
}}}}
\put(5041,-151){\makebox(0,0)[lb]{\smash{{\SetFigFontNFSS{8}{9.6}{\rmdefault}{\mddefault}{\updefault}{$2$}%
}}}}
\put(6751,-151){\makebox(0,0)[lb]{\smash{{\SetFigFontNFSS{8}{9.6}{\rmdefault}{\mddefault}{\updefault}{$2$}%
}}}}
\put(8731,-151){\makebox(0,0)[lb]{\smash{{\SetFigFontNFSS{8}{9.6}{\rmdefault}{\mddefault}{\updefault}{$1$}%
}}}}
\put(9136,-151){\makebox(0,0)[lb]{\smash{{\SetFigFontNFSS{8}{9.6}{\rmdefault}{\mddefault}{\updefault}{$1$}%
}}}}
\put(7111,-151){\makebox(0,0)[lb]{\smash{{\SetFigFontNFSS{8}{9.6}{\rmdefault}{\mddefault}{\updefault}{$1$}%
}}}}
\put(9811,-916){\makebox(0,0)[lb]{\smash{{\SetFigFontNFSS{8}{9.6}{\rmdefault}{\mddefault}{\updefault}{$1$}%
}}}}
\put(7111,-1636){\makebox(0,0)[lb]{\smash{{\SetFigFontNFSS{8}{9.6}{\rmdefault}{\mddefault}{\updefault}{$1$}%
}}}}
\put(5131,-1636){\makebox(0,0)[lb]{\smash{{\SetFigFontNFSS{8}{9.6}{\rmdefault}{\mddefault}{\updefault}{$2$}%
}}}}
\put(3736,-151){\makebox(0,0)[lb]{\smash{{\SetFigFontNFSS{8}{9.6}{\rmdefault}{\mddefault}{\updefault}{$2$}%
}}}}
\put(4636,-1636){\makebox(0,0)[lb]{\smash{{\SetFigFontNFSS{8}{9.6}{\rmdefault}{\mddefault}{\updefault}{$2$}%
}}}}
\put(3691,-1051){\makebox(0,0)[lb]{\smash{{\SetFigFontNFSS{8}{9.6}{\rmdefault}{\mddefault}{\updefault}{$1$}%
}}}}
\put(5716,-826){\makebox(0,0)[lb]{\smash{{\SetFigFontNFSS{8}{9.6}{\rmdefault}{\mddefault}{\updefault}{$1$}%
}}}}
\put(4321,-2671){\makebox(0,0)[lb]{\smash{{\SetFigFontNFSS{8}{9.6}{\rmdefault}{\mddefault}{\updefault}{$1$}%
}}}}
\put(3601,-2356){\makebox(0,0)[lb]{\smash{{\SetFigFontNFSS{8}{9.6}{\rmdefault}{\mddefault}{\updefault}{$2$}%
}}}}
\put(3556,-2671){\makebox(0,0)[lb]{\smash{{\SetFigFontNFSS{8}{9.6}{\rmdefault}{\mddefault}{\updefault}{$4$}%
}}}}
\put(3916,-2671){\makebox(0,0)[lb]{\smash{{\SetFigFontNFSS{8}{9.6}{\rmdefault}{\mddefault}{\updefault}{$2$}%
}}}}
\put(5806,-2896){\makebox(0,0)[lb]{\smash{{\SetFigFontNFSS{8}{9.6}{\rmdefault}{\mddefault}{\updefault}{$2$}%
}}}}
\put(5941,-3391){\makebox(0,0)[lb]{\smash{{\SetFigFontNFSS{8}{9.6}{\rmdefault}{\mddefault}{\updefault}{$1$}%
}}}}
\put(8146,-2851){\makebox(0,0)[lb]{\smash{{\SetFigFontNFSS{8}{9.6}{\rmdefault}{\mddefault}{\updefault}{$1$}%
}}}}
\put(7021,-3436){\makebox(0,0)[lb]{\smash{{\SetFigFontNFSS{8}{9.6}{\rmdefault}{\mddefault}{\updefault}{$1$}%
}}}}
\put(7021,-2221){\makebox(0,0)[lb]{\smash{{\SetFigFontNFSS{8}{9.6}{\rmdefault}{\mddefault}{\updefault}{$1$}%
}}}}
\put(8821,-2761){\makebox(0,0)[lb]{\smash{{\SetFigFontNFSS{8}{9.6}{\rmdefault}{\mddefault}{\updefault}{$1$}%
}}}}
\put(8821,-2221){\makebox(0,0)[lb]{\smash{{\SetFigFontNFSS{8}{9.6}{\rmdefault}{\mddefault}{\updefault}{$2$}%
}}}}
\put(9901,-3346){\makebox(0,0)[lb]{\smash{{\SetFigFontNFSS{8}{9.6}{\rmdefault}{\mddefault}{\updefault}{$1$}%
}}}}
\put(9901,-2761){\makebox(0,0)[lb]{\smash{{\SetFigFontNFSS{8}{9.6}{\rmdefault}{\mddefault}{\updefault}{$1$}%
}}}}
\put(9901,-2221){\makebox(0,0)[lb]{\smash{{\SetFigFontNFSS{8}{9.6}{\rmdefault}{\mddefault}{\updefault}{$2$}%
}}}}
\put(10801,-2221){\makebox(0,0)[lb]{\smash{{\SetFigFontNFSS{8}{9.6}{\rmdefault}{\mddefault}{\updefault}{$1$}%
}}}}
\put(10801,-2761){\makebox(0,0)[lb]{\smash{{\SetFigFontNFSS{8}{9.6}{\rmdefault}{\mddefault}{\updefault}{$1$}%
}}}}
\put(5491,-2221){\makebox(0,0)[lb]{\smash{{\SetFigFontNFSS{8}{9.6}{\rmdefault}{\mddefault}{\updefault}{$1$}%
}}}}
\put(6706,-1636){\makebox(0,0)[lb]{\smash{{\SetFigFontNFSS{8}{9.6}{\rmdefault}{\mddefault}{\updefault}{$2$}%
}}}}
\put(8731,-1636){\makebox(0,0)[lb]{\smash{{\SetFigFontNFSS{8}{9.6}{\rmdefault}{\mddefault}{\updefault}{$1$}%
}}}}
\put(10801,-151){\makebox(0,0)[lb]{\smash{{\SetFigFontNFSS{8}{9.6}{\rmdefault}{\mddefault}{\updefault}{$1$}%
}}}}
\put(10801,-1636){\makebox(0,0)[lb]{\smash{{\SetFigFontNFSS{8}{9.6}{\rmdefault}{\mddefault}{\updefault}{$1$}%
}}}}
\end{picture}%

%% file: leg.pstex_t
\begin{picture}(0,0)%
\includegraphics{leg.pstex}%
\end{picture}%
\setlength{\unitlength}{4144sp}%
\begingroup\makeatletter\ifx\SetFigFontNFSS\undefined%
\gdef\SetFigFontNFSS#1#2#3#4#5{%
  \reset@font\fontsize{#1}{#2pt}%
  \fontfamily{#3}\fontseries{#4}\fontshape{#5}%
  \selectfont}%
\fi\endgroup%
\begin{picture}(3722,603)(1282,68)
\put(1303,132){\makebox(0,0)[lb]{\smash{{\SetFigFontNFSS{12}{14.4}{\rmdefault}{\mddefault}{\updefault}{$0$}%
}}}}
\put(2203,132){\makebox(0,0)[lb]{\smash{{\SetFigFontNFSS{12}{14.4}{\rmdefault}{\mddefault}{\updefault}{$1$}%
}}}}
\put(3103,132){\makebox(0,0)[lb]{\smash{{\SetFigFontNFSS{12}{14.4}{\rmdefault}{\mddefault}{\updefault}{$2$}%
}}}}
\put(1297,512){\makebox(0,0)[lb]{\smash{{\SetFigFontNFSS{12}{14.4}{\rmdefault}{\mddefault}{\updefault}{$n$}%
}}}}
\put(4903,132){\makebox(0,0)[lb]{\smash{{\SetFigFontNFSS{12}{14.4}{\rmdefault}{\mddefault}{\updefault}{$l$}%
}}}}
\put(4897,512){\makebox(0,0)[lb]{\smash{{\SetFigFontNFSS{12}{14.4}{\rmdefault}{\mddefault}{\updefault}{$d_l$}%
}}}}
\put(3097,512){\makebox(0,0)[lb]{\smash{{\SetFigFontNFSS{12}{14.4}{\rmdefault}{\mddefault}{\updefault}{$d_2$}%
}}}}
\put(2197,512){\makebox(0,0)[lb]{\smash{{\SetFigFontNFSS{12}{14.4}{\rmdefault}{\mddefault}{\updefault}{$d_1$}%
}}}}
\end{picture}%

%% file: splayedleg.pstex_t
\begin{picture}(0,0)%
\includegraphics{splayedleg.pstex}%
\end{picture}%
\setlength{\unitlength}{4144sp}%
\begingroup\makeatletter\ifx\SetFigFontNFSS\undefined%
\gdef\SetFigFontNFSS#1#2#3#4#5{%
  \reset@font\fontsize{#1}{#2pt}%
  \fontfamily{#3}\fontseries{#4}\fontshape{#5}%
  \selectfont}%
\fi\endgroup%
\begin{picture}(4523,1597)(1381,-1358)
\put(2116,-1051){\makebox(0,0)[lb]{\smash{{\SetFigFontNFSS{12}{14.4}{\rmdefault}{\mddefault}{\updefault}{$n_j$}%
}}}}
\put(1396,-556){\makebox(0,0)[lb]{\smash{{\SetFigFontNFSS{12}{14.4}{\rmdefault}{\mddefault}{\updefault}{$\#J$}%
}}}}
\put(3097,-388){\makebox(0,0)[lb]{\smash{{\SetFigFontNFSS{12}{14.4}{\rmdefault}{\mddefault}{\updefault}{$d_1$}%
}}}}
\put(3997,-388){\makebox(0,0)[lb]{\smash{{\SetFigFontNFSS{12}{14.4}{\rmdefault}{\mddefault}{\updefault}{$d_2$}%
}}}}
\put(5797,-388){\makebox(0,0)[lb]{\smash{{\SetFigFontNFSS{12}{14.4}{\rmdefault}{\mddefault}{\updefault}{$d_l$}%
}}}}
\end{picture}%